\documentstyle[12pt]{article}

 \input{amssym}

\begin{document}

\vspace{0,3cm }

\centerline{ \large \bf STEENROD OPERATIONS  AND   HOCHSCHILD HOMOLOGY.}

\vspace{3mm}

\centerline{ \small by}
\vspace{3mm}

\centerline{ {\bf {\sc Bitjong Ndombol and Jean-Claude Thomas}}}

\vspace{3mm}
 \noindent University of  {\sc Dschang}  - Cameroun and  University of {\sc
Angers - CNRS-6093}

 D\'epartement de
Math\'ematiques -  2, Bd Lavoisier
 49045 Angers - France.

\vspace{ 0,3 cm}

\rule{16 cm}{0,2mm}
\vspace{ 0,2 cm}

       \noindent {\bf Abstract.} {\small
       Let $X$ be a simply connected space and ${\Bbb F}_p $ be a prime field.
The algebra of
normalized singular cochains $N^*(X; {\Bbb F}_p)$ admits a natural homotopy
structure  which
induces  natural Steenrod operations  on  the Hochschild
       homology $HH_* N^*(X;{\Bbb F}_p) $ of the space $X$. The primary
purpose of this paper is to prove that    the J.
Jones  isomorphism  $HH_*N^*(X;{\Bbb F}_p)  \cong H ^*(X^{S^1};{\Bbb F}_p) $
identifies theses Stenrood operations with those  defined on the cohomology of
the free loop space with coefficients in
${\Bbb F}_p$. The other goal of this paper is to describe a theoritic model
which allows to do some computations. }

\vspace{ 0,2 cm}
\rule{16 cm}{0,2mm}


\vspace{ 0,2 cm}
\rule{16 cm}{0,2mm}


\vspace{ 0,2 cm}

\rule{16 cm}{0,2mm}
\vspace{ 0,6 cm}

       \noindent {\bf Introduction.}

 Let $A =\{A^i\}_{i\geq 0}$ be a  augmented differential graded algebra
over the field ${\Bbb F}_P$. The homology of the
normalized Hochschild {\bf chain} complex  ${\frak C}_*A$   is called  the
Hochschild homology (with coefficients in $A$) of
$(A,d_A)$ and is  denoted by  $HH_*A$.  Let $X$ be a simply connected
space.  In 1987, J.D.S. Jones \cite{[Jo]} constructed
an  isomorphism of graded vector spaces $ HH_*N^*X \cong H^{*}(X^{S^1};
{\Bbb F}_p   )
\,, $
where $X ^{S ^1}$ denotes the free loop space and  $N^*X $ denotes singular
cochains with coefficients in ${\Bbb F}_p$. In view of this result, one may
ask  the following question:

{\it    Does there exist a subcategory ${\bf DA }'$ of the
category ${\bf DA }$ of differential graded algebras such that:

\noindent (Q.1) for any connected space $X$,   $N ^*X$  is an object of
${\bf DA
}'$,

\noindent (Q.2)  if $A$ is an object of ${\bf DA
}'$ then $HH_*A $ is an ${\cal B}_p
$-unstable algebra,

\noindent (Q.3) the Jones'isomorphism is an homomorphism of  ${\cal A}_p
$-unstable algebras? }

Here ${\cal  B}_p $ denotes the large  mod $p$-Steenrod algebra and
$ {\cal  A}_p ={\cal  B}_p/(P^0=id) $ denotes the  usual
mod
$p$-Steenrod algebra.

\vspace{2mm}
In this paper we give an affirmative answer to (Q.1) and (Q.3) and a partial
answer to  (Q.2). For this purpose we introduce
the  notion of $\pi$-$shc$ algebra where $\pi$ denotes  the cyclic group
of order a fixed prime $p$ (1.6). More precisely we prove:

\vspace{2mm}

{\bf Theorem A.} {\it Let $p$ be a prime and $A$ be the mod $p$ reduction of a
  graded module over
${\Bbb Z}$ and denote by
$\beta$ the Bockstein homomorphism.
If
$A$  is a
$\pi$-$shc$  algebra over ${\Bbb F}_p$ then  for any
$i
\in {\Bbb Z}$ and any
$x \in (HH_* A) ^n $
  there exist  well defined homology classes $
Sq^i(x) \in
(HH_*A) ^{n+i} $ if $ p =2  $ and
$P^i (x) \in (HH_*A) ^{n+2i(p-1)}  $if $ p >2$
such that: $Sq^1( x)= P^i (1_{HA}) =0\,, i \neq 0$ and \\
 $Sq^i(x) =  \left\{ \begin{array}{lll} 0  &\mbox{ if } i>
n\\
  x^2   &\mbox{ if }i=n \end{array} \right.$ and
$\beta ^\epsilon  P^i(x) = \left\{ \begin{array}{lll} 0  &\mbox{ if }
2i+\epsilon > n \,, \epsilon = 0,1\\
 x^p    &\mbox{ if }  n=2i \,, \epsilon =0 \end{array} \right. $ \\
Moreover these classes are natural with respect to homomorphisms of
$\pi$-$shc$
 algebras and satisfy the Cartan formula.}

\vspace{3mm}
\noindent{\bf Theorem B.} {\it

\noindent 1) The algebra $N^*X$ is naturally a $\pi$-$shc$ algebra.  Moreover,
any two natural structural maps defining a $\pi$-$shc$
structure on
$N^*X$ are $\pi$-homotopic.

\noindent 2)  If $X$ is 1-connected then, the Jones'quasi-isomorphism  $
{\frak
C}N^*X \to C^*(X ^{S^1})$ identifies the
algebraic Steenrod operations, defined on
$HH_*N^*X$  by theorem A, with the topological Steenrod operations on
$H^*(X^{S^1}; {\Bbb F}_p)$.}

\vspace{2mm}
In  $\S 3 $   we develop the construction of a convenient model for
a $\pi$-$shc$  algebra.
But, at it can be easily imagined, explicit computations with this model
are not very tractable in general (cf. Example 3.9).

\vspace{2mm}
These two theorems are a generalization of theorem A and B, in
\cite{[BT]}. In fact, we proved  that
the Jones' homomorphism  induces a commutative diagram

\centerline{$
\begin{array}{ccccc}
N^* X  &\stackrel{\iota}{\longrightarrow}& {\frak C}N^*X
&\stackrel{\rho}{\longrightarrow}& BN^*X\\
||  && J_X\downarrow && \qquad \downarrow \overline J_X \\
N^*X & \stackrel{\iota '}{\longrightarrow}& N^* X ^{S^1} & \stackrel{\rho'
}{\longrightarrow}& N^* \Omega X
\end{array}
$}

\noindent
where the horizontal arrows in the lower lines are induced from the canonical
fibration
$\Omega X \to X ^{S^1} \to X$ and
$ HJ_X$  (resp. $H\overline J_X$) is an isomorphism of commutative graded
algebras (resp. of commutative graded Hopf
algebras). Now following the lines  of the proof of theorem B  and of
\cite{[BT]}-II-$\S 3$
we  obtain:

\vspace{3mm}

\noindent{\bf Theorem C.} {\it If $X$ is 1-connected then  the
homomorphisms, $H\iota$, $H\rho $,  \ $HJ_X$ and $H\overline
J_X$  respect Steenrod operations.}

\vspace{2mm}

The paper will be  organized as follows:

$\S$ 1 - Proof of theorem A.

$\S$ 2 - Prof of theorem B

$\S$ 3 - $\pi $ shc model.

Appendix A - Technicalities on $shc$ algebras.

Appendix B - Equivariant acyclic model theorem.

\vspace{3mm}

We would like to thank  Katsuhiko Kuribayashi and   Luc Menichi for
helpfull conversations during the preparation of this paper.

\vspace{1cm }

\centerline{ \bf \large $\S$1  Proof of theorem A.}

\vspace{3mm }

\noindent{\bf 1.1 Preliminaries.}       Throughout the paper, $p$ is a
fixed prime, and we work
over a field  ${\Bbb F}_p$. If  $\pi  $ is any finite group,
 the group ring  ${\Bbb F}_p [\pi] $  is a augmented algebra. If two
$\pi$-linear maps $f$ and $g$
are
$\pi$-linear homotopic we write $f\simeq_{\pi}g$.

\vspace{3mm }

\noindent{\bf 1.2. Algebraic Steenrood operations.}  The material involved
in  this section is
contained  in   \cite{[Ma]}. Let  $\pi = \{1,\tau, ...,\tau ^{p-1} \} $ be
the cyclic
group of order $p$.  We identify $\tau$
with the $p$-cycle $(p,1,2,...,p-1)$ of ${\frak
S}_p$ thus $\pi$ acts on $A^{\otimes p}$.  Let $W \stackrel{\epsilon
_W}{\to } {\Bbb F}_p $ be a
projective resolution of ${\Bbb F}_p $ over
${\Bbb F}_p [\pi ]$: $ W = \{W_i\}_{i \geq 0}\,, \quad  \partial : W_i \to
W _{i-1}  \,, W_0\cong
{\Bbb F}_p[\pi ]$
where each $W_i$ is a right projective $ \pi $-module and $\partial $ is $\pi
$-linear. We choose a linear map  $\eta : {\Bbb F}_p\to W $ such that
$\epsilon_W \circ \eta
= id_{{\Bbb F}_p}$. Necessarily, such $\eta $
satisfies also $ \eta \circ \epsilon_W \simeq id _W$.

Let $A =\{ A_i \}_{i \in {\Bbb Z}} $ be a differential  graded algebra (not
necessarily associative). We denote by  $
m^{(p)}$ (resp. $ (Hm)^{(p)}$) the iterated product  $a_1\otimes a_2\otimes
...\otimes  a_p \mapsto a_1(a_2(...a_p))...) $
(resp.  the  iterated product induced on $HA$  by $ m^{(p)} $). Assume that
$\pi$ acts  trivially on $A$ and diagonally on $W \otimes A^{\otimes p}$.
If   there exists   a $\pi$-chain map  $\theta $
such that the left hand diagram induces the right hand   diagram which is
commutative:

\centerline{
$
\begin{array}{lcl}
W \otimes A^{\otimes p} &\hspace{-8mm} \stackrel {\theta}{\longrightarrow}
&\hspace{-8mm} A \\
\hspace{-8mm}\eta  _W \otimes id \uparrow  &\hspace{8mm} \nearrow  m ^{(p)} \\
\hspace{8mm} A ^{\otimes p}
\end{array}
 \qquad
\begin{array}{lcl}
H(W \otimes A ^{\otimes p}) & \hspace{-15mm}\stackrel
{H\theta}{\longrightarrow}   & \hspace{-18mm}HA \\
\qquad \cong \uparrow &   \hspace{6mm}\nearrow \hspace{4mm}(Hm) ^{(p)} \\
\hspace{8mm} (HA)^{\otimes p}
\end{array}
$}

\noindent then  for any
$i
\in {\Bbb Z}$ and any
$x \in H^n A $  there exist  well defined homology classes

\centerline{$
P^i(x) \in
\left\{ \begin{array}{lll} &H ^{n+i}A  &\mbox{ if } p =2 \\
 &H ^{n+2i(p-1)}A  &\mbox{ if } p >2
\end{array} \right.
\,, \quad  \mbox{ and if } p >2 \,, \quad    \overline{P}^i(x) \in H
^{n+2i(p-1)+1 }A\,, $}

\noindent such that:

1) $P^i (1_{HA}) =0\,, i \neq 0$,

2) if $p=2$,  $P^i(x) = 0 $ if $i > n $  and $P^n (x) = x ^2$,

3) if $p>2$,  $P^i(x) = 0 $  if $2i > n $   and
$P^i (x) = x ^p$ if $n=2i$, and
$\overline{P}^i(x) = 0 $  if $2i \geq  n $,

4) If $A$ is the mod $p$ reduction of a free graded module over ${\Bbb Z}$
then
$\beta
\circ P^{i-1}   =i P^{i} $ if $p=2$ and if $p>2$ $
\overline{P}^i =
\beta \circ P ^i $.\\
 Moreover these classes do not depend on the choice of  $W$ neither on
$\eta$, and    are natural with respect to
homomorphisms of  algebras commuting with the structural maps  $\theta$.
These  {\it algebraic Steenrod operations}
\underline{do not} in general satisfy:  $P^i(x)=0$ if $i<0$, $P^0(x)=x $,
the Cartan formula and the Adem relations.
\vspace{3mm}

\noindent{\bf 1.3 Cartan formula.} Assume $A^ i =0\,, i<0$. The structural
map $\theta
: W
\otimes A^{\otimes p}
\to A
$ induces
$
\tilde \theta : A^{\otimes p} \to \mbox{Hom}(W, A)  \,, \quad \tilde \theta
(u)(w)= (-1)
^{|u|\, |w|} \theta (w\otimes u)\,, \quad u \in
a^{\otimes p}\,, w \in W \,.
$
Let us precise that if $f\in \mbox{Hom} ^k (W,A)= \displaystyle \prod
_{i\geq 0} \mbox{Hom}
(W_i, A^{k-i}) = \displaystyle \bigoplus
_{i=0}^k  \mbox{Hom} (W_i, A^{k-i}) $ then: $
Df = d\circ f - (-1)^k f \circ \partial$, $  (\sigma f)(w) = f(w\sigma )\,,
\quad \sigma \in \pi\,, w \in W\,.
$ A direct verification shows that $  \mbox{Hom} ^k (W,A)$ is a
$\pi$-complex and that
$\tilde \theta $ is a $\pi$-chain map. If $ ev_0 : \mbox{Hom}(W, A) \to A$
denotes the evaluation map on the  generator
$e_0= \eta (1) $ of $W_0$,  the left hand diagram  induces the right hand
diagram which is commutative.

\centerline{$
\begin{array}{ccc}
A^{\otimes p} & \stackrel {\tilde \theta }{\longrightarrow}& \mbox{Hom} (W,
A)\\
\hspace{5mm} m ^{(p)}& \searrow & \hspace{8mm} \downarrow ev _{0}\\
&&A
\end{array}
 \quad
\begin{array}{ccc}
(HA)^{\otimes p} & \stackrel {H\tilde \theta }{\longrightarrow}&H \left(
\mbox{Hom} (W,
A)\right)\\
\hspace{8mm} (Hm)^{(p)}&\searrow  & \hspace{6mm}\downarrow H(ev _{0})\\
&& HA
\end{array}\,.
$}

\noindent Let
$
\psi : W
\to W
\otimes W$  be  a  diagonal approximation and let $m$
denotes  the product in the  algebra $A$. The formula $f\cup g = m \circ
(f\otimes g) \circ \psi$ defines  cup product
$ \mbox{Hom} ^k (W,A) \otimes  \mbox{Hom} ^l (W,A) \to  \mbox{Hom} ^{k+l}
(W,A)\,, \quad f\otimes g
\mapsto f \cup g
$ and  $\mbox{Hom} (W, A)$ is   (non associative)  differential graded algebra.

\vspace{3 mm}
\noindent{\bf 1.4  Proposition.} {\it Let $(A,\tilde \theta )$ be as above.
If we assume that $H^*\tilde \theta $ respects
products  then the algebraic Steenrod operations defined by $\tilde \theta
$ satisfy the
Cartan
formula: $
P^i(xy)= \sum _{j+k=i} P^j (x) P^k(y) \,, \quad x,y \in H^*A\,. $}

\vspace{2 mm}

\noindent {\bf Proof.} We consider the {\it standard small free resolution}
of a finite cyclic group
\cite{[Eil]}: the  (right) $\pi$-free acyclic chain
complex defined by:
$ W= \{W_i\}_{i \geq 0} $, $  W_i = e_i {\Bbb F}_p [\pi]$, $\partial : W_i
\to W_{i-1} $ and if
$\tau $ is a fixed generator of $ \pi $ then
$
  \partial e_{2i+1} = (1+\tau) e_{2i}$, $  \partial e_{2i}= (1+\tau+...+ \tau
^{p-1}) e_{2i-1} $. The augmentation
$\epsilon _W : W \to {\Bbb F}_p $ is defined by $  \epsilon_W (e_i)= 0\,,
i>0 $ and $
\epsilon_W (e_0)=1 $.  Let $\triangle A \subset A ^{\otimes p}$ be the
subspace generated by the  image of the diagonal map
 $\triangle: A \rightarrow A^{\otimes p} $
 $a \mapsto  a^{\otimes p} $. Thus  $\pi$ acts trivially on the sub-algebra
 $\triangle A $, and there exists a $\pi$-module $M\subset A ^{\otimes p}$
such that,  \cite{[Ma]}-Lemma 1.3,
$
H(W\otimes _{\pi} A ^{\otimes p})= \left( \bigoplus _{i= 0} ^\infty {\Bbb
F}_p e_i
 \otimes \triangle A \right) \oplus \left(
{\Bbb F}_p e_0  \otimes M \right) \,.
$
The  algebraic
Steenrod operations are defined by  (\cite{[Ma]}-Proposition 2.3): for  $x
\in H^n A$  if $p=2$ by,  $P^i(x)= \theta ^*
(e_{n-i}
\otimes x ^{\otimes p}):= Sq^i(x)$ and if $p$ is odd by ,
$P^i(x)=
(-1) ^i \nu (n) \theta ^* (e_{(n-2i)(p-1) } \otimes x ^{\otimes p})$, $
\overline{P}^i(x)=
(-1) ^i \nu (n) \theta ^*  (e_{(n-2i)(p-1)-1 } \otimes x ^{\otimes p}) $.
Here $\nu (n) = (-1)^j \left( \left(
\frac{p-1}{2}\right) !\right)^\epsilon  $  if $n=2j+\epsilon $, $\epsilon =
0,1$ and $\theta ^*$ denotes the map induced by
$\theta $ on $ H(W\otimes _{\pi} A ^{\otimes p})$.  A  diagonal approximation
is explicitely given by, \cite { [C-E]} -$\S 31$-Chap. $X\!I\!I$:  $\psi
e_i (\tau ^s )= \psi(e_i) \tau ^s $,
$
\psi (e_{2l+1} ) = \sum _{j+k=l} e_{2j} \otimes e_{2k+1} +  \sum _{j+k=l}
e_{2j+1} \otimes
e_{2k}\tau ^{-1}$ and $ \psi (e_{2l} ) =
\sum _{j+k=l} e_{2j} \otimes e_{2k} +  \sum _{j+k=l-1} \sum _{0\leq r<s<p}
e_{2j+1} \tau ^{-r} \otimes e_{2k+1 }\tau ^{-s}$.  For simplicity, let us
assume
that
$p=2$. First observe that the restriction of  $\tilde \theta $ on
$\triangle A \subset A ^{\otimes 2}$
 factorizes through a map $\tilde \theta_{\pi}:
\triangle A  \rightarrow \mbox{Hom}_{\pi}(W, A) $ where
$\mbox{Hom}_{\pi}(W, A) $ denotes the sub-complex
of
$\pi$-linear maps. Secondly, note  that $H^*(\mbox{Hom}_{\pi}(W, A))\cong
H^*(\mbox{Hom}_{\pi}(W,H^* A))$
and that the map $H^*\tilde \theta $ factorizes through a map
$(H^*\tilde \theta )_{\pi}:
K\cap kerd \rightarrow \mbox{Hom}_{\pi}(W, H^*A) $. Let us denote by $j:
\mbox{Hom}_{\pi}(W, A)\rightarrow \mbox{Hom}(W, A)$
the natural inclusion  of complexes.  By definition, if $x \in H ^n A $ and
$a$ denotes  any
cocycle such that $x=cl(a)$: $
Sq^ix = H^*j\circ ((H^*\tilde \theta )_{\pi}(a^{\otimes 2}))(e_{n-i}) =cl (
\tilde \theta ( a^{\otimes 2})(e_{n-i}))\,.
$
 If   $y=cl(b) \in H^m A$ then:\\
$
\begin{array}{lll}
Sq ^i (xy)&=
cl\left( \tilde \theta (ab \otimes ab ) (e_{n+m-i})\right)
= cl\left( \tilde  \theta ((a\otimes a)(b \otimes b) )
(e_{n+m-i}) \right)
\\
 &
= cl\left(  (\tilde \theta (a\otimes a) \cup \tilde \theta (b \otimes b) )
(e_{n+m-i})\right)
= cl \left( \sum _{j=0}^{n+m-i} \tilde \theta (a\otimes a)  (e_{j})
\tilde \theta (b\otimes b)  (e_{m+n-i -j} \tau ^j)\right) \\
&= \sum _{j=0}^{n+m-i}  cl \left( \tilde \theta (a\otimes a)  (e_{j})
\right)  cl \left(
\tilde \theta (b\otimes b)  (e_{m+n-i -j} )\right) = \sum _{j=0}^{n+m-i}
Sq^{n-j} x Sq ^{-n+i+k} y
\\
 &= \sum _{j+k =i} Sq^{j} x Sq^{k} y\,.
\end{array}
$

\vspace{3mm}

\noindent{\bf 1.5. The differential graded algebra ${\frak C}_*(A)$.} The
bar and the cobar construction can be regarded as a
pair of adjoint functors
$
B : \mbox{ \bf DA}  \leftrightarrow \mbox{ $\cal N$-{\bf DC}} : \Omega \,,
$
where { \bf DA} is the category of suplemented differential graded
\underline{associative}  algebras and
  $\cal N $-{\bf DC} is the category of
locally conilpotent supplemented graded  \underline{coassociative}
coalgebras. (See for instance,
\cite{[FHT]}-2.14 for more
details.) This leads to a  natural homomorphism
$
\alpha _A : \Omega BA \to A
$
of differential graded algebras.   The elements of $BA $ (resp. $\Omega C$)
are denoted
$
[a_1|a_2|...|a_k] \in B_kA\,,$ and $[\,]:=1 \in B_0A \cong  {\Bbb F}_p$
(resp. $<c_1|c_2|...|c_k> \in \Omega_k A\,,$ and
$<\,>:=1
\in \Omega _0 A \cong  {\Bbb F}_p$). We remind the reader that the linear map
$
i_A : A={\Bbb F}_p\oplus \overline{A} \to \Omega BA \,, \quad i_A(1)=1
\mbox{ and } i_A(a)= s^{-1}s a=<[a]>
\,,  a \in \overline A \,,
$
commutes with the differentials and satisfies: $ \alpha _A \circ i_A = id $,
$ i_A \circ \alpha _A \simeq  id_{\Omega B A}$.

 Recall also  that  the {\it Hochschild complex} of $A$, is the graded
vector space  $
 {\frak C}_*A = \{ {\frak C}_kA\}_{k \geq 0} \,, \quad  {\frak C}_kA   = A
\otimes B_k A\,,
$ where  a generator of  ${\frak C}_kA$ is of the form
$a_0[a_1|a_2|...|a_k]$ if $k>0$
and  $a[\,]$ if $k=0$. We set $\epsilon _i = |a_0|+
|sa_1|+|sa_2|+...+|sa_{i-1}|$, $i\geq 1$
 and define the differential $d=d^1+d^2$ by:

\centerline{$
\begin{array}{ll}
d^1 a_0[a_1|a_2|...|a_k] &= da_0[a_1|a_2|...|a_k] - \sum _{i=1}^k
(-1)^{\epsilon
_i} a_0[a_1|...|da_i|...|a_k]\\
d^2 a_0[a_1|a_2|...|a_k] &=  (-1) ^{|a_0|} a_0a_1[a_2|...|a_k] +
 \sum _{i=2}^k (-1)^{\epsilon _i} a_0[a_1|...|a_{i-1}a_i|...|a_k]\\
& \qquad -  (-1)^{|sa_k|\epsilon _k} a_ka_0[a_1|...|a_{k-1}]
\end{array}
$}
\noindent  Henceforth we shall discard the lower star in ${\frak C}_*A$ and
write ${\frak
C}A$.
 The associatavity of the algebra  $A$ is necessary to prove that $d\circ d
=0$ in $BA$ as well  as in ${\frak C}A$ and if
$A$ is commutative then $BA$  (resp. ${\frak  C}A$) is a differential
graded Hopf algebra (resp. a differential graded
algebra.)  Since we are interested to the case $A= N ^*(X)$, and since this
differential graded algebra is not commutative but
only {\it strongly homotopy commutative} (a $shc$ algebra  for short),
\cite{[Mu]},   it seems  natural for our purpose, to
consider
$shc$ algebra    in order to define a product on ${\frak  C}A$.

 A  $shc$ {\it algebra }
 is a pair $(A,\mu)$ such that
$
\mu :  \Omega B (A^{ \otimes 2})   {\rightarrow} \Omega B A
$
 is a homomorphism of differential graded algebras where
$ \alpha \circ \mu \circ i _{A \otimes A}= m_A$  ($m_A$  denotes  the
product in  $ A $). Moreover,
$\mu
$ satisfies the unity axiom, the associativity axiom and the commutativity
axiom
as described in \cite{[Mu]}-4.1. Recall that two homomorphisms of
differential graded algebras  $f, g : A \to A' $
{\it are
homotopic}   (we write $f\simeq_{{\bf DA}}g$, even if $A$ and $A'$ are not
associative) if there exists a $(f, g)$-derivation
$h$ such that
$f-g  = dh + hd$. Let
$(A ,d, \mu )$ and
$(A',d',
\mu') $ be two
$shc$  algebras. A  {\it strict  $shc$  homomorphism} $f:   (A ,d, \mu ) \to
(A',d', \mu') $   is a morphism $f \in {\bf DA}(A,A') $ such that  $ \Omega
B f \circ \mu \simeq _{\bf DA} \mu ' \circ
\Omega B(f
\otimes  f)$.

The chain map  ${\frak C}(\alpha
_A)$ is  surjective  and
 $H{\frak C}(\alpha _A)= HH_*(\alpha _A)$ is an isomorphism. We denote by
$ s_A : {\frak C}A \to {\frak C}\Omega BA  $
 any chain  map such that
 ${\frak C}(\alpha _A) \circ s_A = id $ and $ s_A \circ {\frak C}(\alpha _A)
\simeq id$. When $(A,\mu)$ is a $shc$  algebra   a product $m_{{\frak C}A}
:  {\frak C}  A\otimes
       {\frak C }A \rightarrow {\frak C} A$
is defined by  the  composite

\centerline{$
 {\frak C} A \otimes {\frak C} A
       \stackrel{sh}{\rightarrow}  {\frak C}(A \otimes A)
       \stackrel{ s_{A\otimes A}}{\rightarrow} {\frak C}
       \Omega B(A\otimes A)
       \stackrel{{\frak C} \mu_A}{\rightarrow} {\frak C} \Omega B A
\stackrel{{\frak C} \alpha _A}{\rightarrow} {\frak C} A
$}

\noindent
where $sh$ denotes the shufle product, \cite{[BT]}.  More precisely, we
proved:{ \it If $A$ is a $shc$  algebra   then

1) $BA$ is a differential graded Hopf algebra such that $HBA$ is a commutative
graded Hopf algebra,

2) ${\frak C}A$ is a (non associative) differential graded algebra  such
that $HH_*(A):= H {\frak C}A$ is a
commutative  graded algebra,

3) If $A'$ is an other $shc$  algebra   and if $f$ is a homomorphism of
$shc$  algebras then there is a commutative diagram

\centerline{
$
\begin{array}{ccccc}
A &\stackrel{\iota}{\longrightarrow}& {\frak C}A
&\stackrel{\rho}{\longrightarrow}& BA\\
f \downarrow && \hspace{-5mm}{\frak C}f \downarrow && \qquad \downarrow Bf\\
A'& \stackrel{\iota '}{\longrightarrow}& {\frak C}A'& \stackrel{\rho'
}{\longrightarrow}& BA'
\end{array}
$}

\noindent where  the  maps $\iota $, $\iota'$, $\rho $, $\rho '$ and
${\frak C}f$   are homomophisms
of differential graded algebras,  $Bf$  is a homomorphism of differential
graded
Hopf algebras.}

If  $A$ is a   $shc$  algebra  and if $p=2$,  then   both $Sq^n$ and
$ Sq ^{n-1}$ are determined on $H^nA$,  (\cite{[BT]}-6.6   or
\cite{[Mu]}-4.8),  while    $Sq^{n-2}$ is not determined on  $(HH_*A)^n$,
(cf. Example 3.8).   Therefore, we  need  to
enrich the notion of $shc$ algebra.

\vspace{3mm }

\noindent{\bf 1.6  $\pi$-$shc$ algebra.}  Let $(A,d,\mu )$ be a $shc$
algebra.
For any  $n \geq 2$, by Lemma A.3,  there exists a homomorphism of
differential graded algebras, called the {\it $shc$ iterated structural map}
$
\mu ^{(n)} : \Omega B \left( A ^{\otimes n} \right) \to \Omega BA
$
such that:  $ \mu ^{(2)} = \mu$ and  $\alpha _A \circ \mu ^{(n)} \circ i_{A
^{\otimes n}}
\simeq m_A  ^{(n)} $. Let $W, \psi$ as in 1.3.

 A    $shc$  algebra  $(A, \mu ) $ is, a {\it $\pi$-$shc$
algebra} if   there exits a map $ \tilde \kappa_A  : \Omega B (A ^{\otimes
p}) \to \mbox{Hom}( W, A) $ which is both a
  $\pi$-linear map and a
homomorphism of differential graded algebras
such that the next   diagram  commutes, up to a derivation homotopy.

\centerline{$
\begin{array}{ccccc}
\Omega B (A ^{\otimes p})&\stackrel{\tilde \kappa_A }{\longrightarrow} &
\mbox{Hom}(W,A)  \\
\hspace{15mm}{\alpha _A \circ \mu  ^{(p)}}& \hspace{8mm}  \searrow \qquad
& \downarrow ev_{0} \\
&&   \hspace{-7mm}  A
\end{array}\,.
$}

\noindent
Here
$\sigma \in  {\frak S}_p$  acts  on
 $B (A ^{\otimes p})$ by the rule
$
\quad \sigma [x_1|x_2|...|x_k] = [\sigma x_1|\sigma x_2|...|\sigma x_k] \,,
\quad x_i\in A ^{\otimes p}  $ and on  $\Omega B (A ^{\otimes p}) $ by the
rule
$
\sigma \in \pi \,, \quad \sigma <y_1|y_2|...|y_l> = <\sigma y_1|\sigma
y_2|...|\sigma y_l> \,, \quad y_j\in B(A^{\otimes p}) $.

A {\it strict $\pi$-$shc$ homomorphism}  $f: (A, \mu_A,\tilde \kappa_A
)\rightarrow (A', \mu_{A'},\tilde
\kappa_{A'} )$ is a strict $shc$  homomorphism  satisfying
$
 \kappa_{A'} \circ \Omega B (f ^{\otimes p})\simeq_{\pi \mbox{-} {\bf DA}}
\mbox{Hom}(W, f)
\circ\tilde \kappa_A
\,.$

\vspace{3mm}

 If $ (A, \mu,
\tilde
\kappa_A )$ is a $\pi$-$shc$ algebra  then   the
next diagram commutes, up to a $\pi$-linear  homotopy:

\centerline{$
\begin{array}{ccc}
\Omega B \left( A ^{\otimes p}\right)  & \stackrel { \tilde \kappa_A
}{\longrightarrow}& \mbox{Hom} (W, A)\\
\alpha  _{A^{\otimes p}} \downarrow& & \downarrow ev _{0}\\
A^{\otimes p}  & \stackrel {m ^{(p)} }{\longrightarrow}&A
\end{array} \,,
$}

\noindent
and  the
Steenrod operations,
defined  for $ x=cl(a) \in H^n A$  by :
$
Sq^i(x) = cl(\tilde \kappa  (i_{A^{\otimes 2}} (a ^{\otimes 2}))(e_{n-i}))$
if $  p=2$ and by
$P^i(x)=(-1) ^i \nu (n)   cl(\tilde \kappa  (i_{A^{\otimes p}}(a ^{\otimes
p}))
(e_{(n-2i)(p-1)}))$ if    $ p>2$,  satisfy the Cartan formula (see Proof of
1.4  for notations).

\vspace{3mm}

 In order to construct Steenrod operations
on
${\frak C}(A)$ we need to compare $Hom(W, {\frak C}A)$ and ${\frak C}
Hom(W, A)$. But this last expression is not
defined unless $Hom(W, A)$ is associative. If $p\neq 2$, the
diagonal approximation considered in the proof of Proposition 1.4   is not
coassociative.  Thus, for $p$ odd,   we will  consider the {\it  standard
resolution}  of a finite cyclic group $\pi$, \cite{[ML]}, \cite{[B]}.  The
Alexander-Whitney diagonal approximation $\psi_W
$,
admits the counity $\epsilon $ and  is strictly coassociative.

\vspace{3mm}

\noindent{\bf 1.7. Proposition.} {\it If $A$ is a differential graded
algebra and if $(W, \psi _W)$ is a coassaciative
coalgebra   then there
exists a natural  chain map
$\phi_A
$
such that  the following diagram commutes

\centerline{$
\begin{array}{cccccc}
 {\frak C}\mbox{\rm Hom}(W,A)&\stackrel
{\phi_A}{\longrightarrow} & \mbox{\rm Hom} (W, {\frak C} A)\\
{\frak C}ev_{0} \searrow  & &\swarrow ev_{0} \\
& {\frak C} A &
\end{array}
$}

\noindent
Moreover, $\phi $ is $\pi $-linear ($\pi$ acts trivially on $B= A\,, {\frak
C} A $, diagonally on $\mbox{\rm Hom}(W,B)$).}

\vspace{2mm}

\noindent{ \bf Proof.} We define $\phi_A $ by:$ \phi_A (f_0[])= f_0$ and if
$k>0$ by $
\phi_A (f_0[f_1|f_2|...|f_k])=
(id \otimes s^{\otimes k}) \circ (f_0\otimes  f_1\otimes f_2\otimes ...
\otimes f_k) \circ \psi_W ^{(k)}$
where  $\psi_W ^{(k)}$ denotes the iterated diagonal:  $\psi_W
^{(1)} = \psi_W $,  $\psi_W
^{(k+1)}= (id \otimes \psi_W) \circ \psi_W
^{(k)}$ and  $s : A \to sA $ denotes  the suspension of degrees. Obviously,
$\phi_A$ is natural in $A$. Let us check in detail that
$\phi_A $ commutes with the differentials. ( $\sigma_k$  denotes the cycle
$(1,..,k)\to (k,1,..., k-1)$)

\centerline{$
\begin{array}{rlll}
\phi_A (d^1&(f_0[f_1|f_2|...|f_k])) \\
&=
\phi_A \left( df_0[f_1|f_2|...|f_k] - \sum _{i=1} ^k  (-1)^{\epsilon _i}
f_0[f_1|...|df_i|...|f_k] \right)\\
&=
(id \otimes s ^{\otimes k})
\left( df_0\otimes f_1\otimes f_2\otimes ...\otimes f_k - \sum _{i=1} ^k
(-1)^{\epsilon _i} f_0 \otimes
f_1\otimes ...\otimes df_i\otimes ...\otimes f_k
\right) \circ \psi_W ^{(k)} \\
&=
(id \otimes s ^{\otimes k})
\left( (d_A\circ f_0- (-1) ^{|f_0|} f_0 \circ \partial )\otimes f_1\otimes
f_2\otimes ...\otimes f_k \right. \\
& \left. \qquad - \sum _{i=1} ^k
(-1)^{\epsilon _i} f_0
\otimes f_1\otimes ...\otimes (d_A\circ f_i- (-1) ^{|f_i|} f_i \circ \partial
)\otimes ...\otimes f_k
\right) \circ \psi_W ^{(k)} \\
&=
(id \otimes s ^{\otimes k})
\left( (d_A  \otimes id -id \otimes d _{A^{\otimes k}}) \circ   (f_0\otimes
f_1\otimes
f_2\otimes ...\otimes f_k) \circ \psi_W ^{(k)} \right.  \\ &  \left. \qquad
 -
(-1) ^{|f_0] + |f_1|+...+|f_k|+k } (id \otimes s ^{\otimes k})  \circ (
 f_0
\otimes f_1\otimes ...\otimes  f_i\otimes ...\otimes f_k) \circ \partial _{W
^{\otimes k}}
\right) \circ \psi_W ^{(k)} \\
&=
d^1_{{\frak C}A} \circ \phi_A ( (f_0[f_1|f_2|...|f_k]) - (-1) ^{\epsilon _k}
\phi_A ( (f_0[f_1|f_2|...|f_k]) \circ \partial
\end{array}
$}

\centerline{$
\begin{array}{rlll}
\phi_A (d^2(f_0[f_1|f_2|...|f_k]))  &=
\phi_A \left( (-1)^{|f_0|} (f_0\cup f_1)[f_2|...|f_k] + \sum _{i=2} ^k
(-1)^{\epsilon _i} f_0[f_1|...|f_i \cup f_i |...|f_k]
\right.\\
& \qquad \qquad \left. + (-1)^{|sf_k|\epsilon _k } (f_k\cup
f_0)[f_1|...|f_{k-1}]
\right)\\
 &=
(id \otimes s ^{\otimes k-1}) \circ
\left( (-1)^{|f_0|} ( m_A \circ (f_0\otimes f_1) \circ \psi_W ) \circ
(f_2\otimes
...\otimes f_k]) \right.
\end{array}
$}

\centerline{$
\begin{array}{rlll}
& \qquad
 + \sum _{i=2} ^k
(-1)^{\epsilon _i} f_0[f_1|...|m_A \circ (f_i \otimes f_i) \circ \psi_W
|...|f_k]+ \\ & \qquad \qquad \left.
(-1)^{|sf_k|\epsilon _k } (m_A \circ (f_k\otimes  f_0) \circ
\psi_W)[f_1|...|f_{k-1}]
\right) \circ \psi_W ^{(k-1)} \\
 &=
(id \otimes s ^{\otimes k-1}) \circ
\left( ( m_A \otimes id
 + \sum _{i=2} ^k
(-1)^{i} id \otimes m_A \otimes id + (-1)^{k+1} ( m_A \otimes id)\circ
\sigma_k
\right)
 \\ & \qquad \qquad  \qquad  \circ (1\otimes s ^{\otimes k} ) ^{-1} \circ
\phi_A
((f_0[f_1|f_2|...|f_k])) \\
&= d ^2 _{{\frak C}A} \phi_A (f_0[f_1|f_2|...|f_k]) \,.
\end{array}
$}

\noindent
Observe that we need also  the coassociativity property of the coproduct
$\psi_W$  in the proof
of  the second equality. It is straightforward
to  check that the remaining  claims in 1.7  hold.

\vspace{3mm }

\noindent{\bf 1.8  End of the proof of Theorem A.}  One  also define, for $
n \geq 3$,   $(k_1,k_2,...k_n)$-shuffles,
the {\it
$n$-iterated shuffle map} $
sh ^{(n)} : ({\frak C}A)^{\otimes n} \to {\frak C}(A^{\otimes n}) $  and
the  iterated product
$m_{{\frak C}A}^{(p)}  :  {\frak C}  A^{\otimes p} \rightarrow {\frak C} A$:

\centerline{$
 {\frak C} A ^{\otimes p}
       \stackrel{sh^{(p)} }{\rightarrow}  {\frak C}(A^{ \otimes p})
       \stackrel{ s_{A^{ \otimes p}}}{\rightarrow} {\frak C}
       \Omega B(A^{\otimes p})
       \stackrel{{\frak C} \mu^{(p)}_A}{\rightarrow} {\frak C} \Omega B A
\stackrel{{\frak C} \alpha _A}{\rightarrow} {\frak C} A \,,
\qquad\mbox{(compare with 1.5)}.
$}

 Consider $\phi _A $ as defined above, then applying the
functor ${\frak C}$ to the diagram that appear in the definition of a
$\pi$-$shc$ algebra (1.6),   we obtain
the following  diagram:

\centerline{$
\begin{array}{cccccc}
{\frak C} \Omega B (A^{\otimes p})& \stackrel
{ \phi_A  \circ {\frak C}\tilde \kappa_A }
 {\longrightarrow} & \mbox{Hom} (W, {\frak C} A)\\
 \hspace{-12mm} {\frak C}( \alpha _{A^{\otimes p}}) \downarrow &    & \\
{\frak C} (A^{\otimes p})&&\qquad \downarrow ev_{0}&\\
\hspace{-5mm} sh ^{(p)} \uparrow && \\
({\frak C} A)^{\otimes p}  &\stackrel{ m _{{\frak C}A} ^{(p)} }{
\longrightarrow}  & {\frak C} A\,.
\end{array}
$}

\noindent
Since $ev_0\circ \tilde \kappa_A \simeq_{{\bf DA}} \alpha_A\circ
\mu_A^{(p)}$, Lemma A.6
implies  that
 $ev_0\circ \phi_A \circ{\frak C}\tilde \kappa_A \simeq  {\frak C}\mu
^{(p)}$. Now observe that in the sequence $
{\frak C}\Omega
B (A^{\otimes p}) \stackrel{ {\frak C}(\alpha _{A^{\otimes p}})
}{\longrightarrow } {\frak C} (A^{\otimes p})
\stackrel{(sh ^{(p)})   }{\longleftarrow } \left({\frak C} A\right)
^{\otimes p}$
of   chain
maps:

a) ${\frak C}(\alpha _{A^{\otimes p}})$ is $\pi$-linear. Indeed
  $\sigma \in  {\frak S}_p$  acts on ${\frak C} (A ^{\otimes p})$ (resp. on
${\frak
C}\Omega B ( A ^{\otimes p})$   by the rule
$
 \quad \sigma (x_0[x_1|x_2|...|x_k]) = \sigma x_0[\sigma x_1|\sigma
x_2|...|\sigma x_k] \,, \quad
x_i\in   A ^{\otimes p}$ (resp. $x_i \in \Omega B A ^{\otimes p}$). Then,
  ${\frak C}\Omega
B (A^{\otimes p})$ and $  {\frak C} (A^{\otimes p}) $ are ${\frak
S}_p$-complexes, $ \alpha _{ A ^{\otimes p}}$ and thus ${\frak C} \alpha
_{A ^{\otimes p}}$
are
${\frak S}_p$-linear.

b)   $sh ^{(p)}$ is also ${\frak S}_p$-linear. Indeed,  one simply
writes:
$
sh^{(p)} (a^1_0[ a^1_1 \vert a^1 _2\vert ... \vert a^1_{n_1}] \otimes
a^2_0[ a^2_1
\vert ... \vert a^2_{n_2}] \otimes ...
\otimes a^p_0[ a_1  ... \vert a^p _{n_p}] )
 = (a^1 _0\otimes a^2 _0 \otimes ... \otimes  a^p _0 ) \sum_{ \sigma }
\sigma [a^1 _1\otimes 1\otimes 1 ... \otimes 1|... | a^1 _{n_1} \otimes 1
\otimes 1 ... \otimes 1 |
1 \otimes a^2 _1\otimes 1\otimes 1 ... \otimes 1 |...|
1\otimes 1\otimes 1 ... \otimes 1 \otimes a ^p_i |...| 1\otimes 1\otimes 1 ...
\otimes 1 \otimes a ^p_{n_p}] $
where $\sigma $ ranges over all $(n_1,n_2,...,n_p)$-shuffles and  observe
that the set of all
$(n_1,n_2,...,n_p)$-shuffles coincides with the set of all $(\lambda
(n_1),\lambda (n_2),..., \lambda (n_p))$-shuffles
when $\lambda \in {\frak S}_p$.

c)   ${\frak C}(\alpha _{A^{\otimes p}}) $ is a homomorphism of
differential graded algebras. This
 follows directly  from  F4, F6 and  F8.

From, properties a), b), c) and  following the lines of the proof on
\cite{[Ma]}-Proposition 2.3, (sketched in
the proof of Proposition 1.4), using the section $s_{A^{\otimes p}}$, this
diagram  defines natural  Steenrod operations on
$HH_*(A)$.   In order to establish the Cartan
formula it is enough, by Proposition 1.4, to prove that $H^*\tilde \kappa_A $
respects products. This  is    Lemma A.5. Naturality   follows directly
from Lemma A.6.

\vspace{1cm }

\centerline{ \bf \large $\S$2 -  Proof of theorem B.}

\vspace{3mm }

\noindent {\bf 2.1 Twisting cochains.}  In \cite {[Mu]}-4.7, H.J. Munkholm
has
established the existence of a natural transformation
$\mu _X : \Omega B (N^* X \otimes   N^*X) \to \Omega B N^*X  $ such that
$\alpha_{N^*X} \circ \mu_X \circ  i_{N^* X \otimes
N^*X}$  is the usual cup product $N^* X \otimes N^*X \to N^*X$. Moreover,
$(N^*
X , \mu_X)$ is a $shc$-algebra. Let us denote by $\mu _X ^{(p)} : \Omega B
\left( (N^*X) ^{\otimes p}\right) \to
\Omega B N^*X  $
the iterated $shc$ structural map and by
$t _X
$ the  twisting cochain associated to $ \alpha _{N^*X} \circ \mu _X ^{(p)}$: $
t_X \in \mbox{Hom} ^1 \left(B \left((N^*X)^{\otimes p}\right),
 N^*X\right)\,,   t \cup t = Dt\,, \quad  t(x)(e_0) = \mu (
<x> )
\,, x
\in
\overline{BN^*X}\,, t([\,])=0
$ where $D$  denotes the differential in
$ \mbox{Hom}  \left(B \left((N^*X)^{\otimes p}\right),  N^*X \right)$ (See
\cite{[Mu]}-1.8 for more details on twisting
cochains).  We want to construct a twisting cochain
$
t'_X \in \mbox{Hom} ^1 \left(B \left((N^*X)^{\otimes p}\right),
\mbox{Hom}(W,N^*X)\right)\,,
 t' \cup t' = D't' $
such that:   $ev_0 \circ t'_X = t_X$, $t'_X$ is natural in $X$ and  $t'_X$
is $\pi $-linear. (Here $D'$ denotes the
differential in $\mbox{Hom} \left(B \left((N^*X)^{\otimes  p}\right),
\mbox{Hom}(W,N^*X)\right)$.) Such a twisting cochain
$t'_X$ will determine  a homomorphism of differential  graded algebras
$
\tilde \kappa _X : \Omega B ( N^*X)^{\otimes p}  \to  \mbox{Hom} (W,N^*X)
$
such that  $( N^*X, \mu_X , \tilde  \kappa _X)$ is a $\pi$-$shc$  algebra.
Moreover, by Corollary B.4,
 $\tilde  \kappa _X$
is well determined up to a $\pi$-linear  homotopy, and so the first part of
theorem B will be proved.

\vskip 3mm

\noindent {\bf 2.2} { \bf Proof of part 1.}   We construct $t'$ by
induction on the degree in $W=
\{W_i\}_{i \geq 0}$ and on the bar degree:
$
\overline{B\left( (N^*X)^{\otimes p}\right) } =\left\{B_k\left((N^*X)^{\otimes
p}\right)\right\}_{k \geq 1}\,.
$   We denote by
$\Delta $ the coproduct in $B
((N^*X) ^{\otimes p})$. For simplicity, if $x \in B_k(N^*X^{\otimes p})$ we
write
$\Delta x = 1\otimes x + x \otimes 1+ \overline \Delta x =  x'\otimes x''$
where
  $\overline \Delta x \in
\oplus _{l'+l''= k}  B_{l'}\left((N^*X)^{\otimes p}\right)\otimes
B_{l''}\left((N^*X)^{\otimes p}\right) $  with  $ l'$
(resp
$l''$)  such that $ 0< l' <k$ (resp. $ 0< l'' <k$).

First, if $e_0$ denotes a generator of $W_0$, we define $t'$ on $W_0
\otimes\overline{B\left( (N^*X)^{\otimes p}\right)
 }$  by setting:

\centerline{$
(\star _0) \qquad t'(x) (e_0 \sigma ) = t(\sigma ^{-1} x) \,,\qquad \sigma
\in \pi\,,
 \quad  x \in B_k N^*X \,. $}

\noindent
$t'$ is a $\pi$-linear map, that is:  if $\sigma \in \pi $ and $  x \in
\overline{B( N^*X)^{\otimes p}}$ then
 $ t'(\sigma x)= \sigma (t'( x )) $. It suffices to prove this last formula
for
the elements  $ e_o \tau ^j$ in $W_0 $:

\centerline{$
\begin{array}{ll}
(\sigma t'(x))(e_0 \tau ^j) &= t'(x)( e_0 \tau ^j \sigma ^{-1})=
t'(x)\left( e_0  \left ( \sigma (\tau ^j)^{-1} \right)^{-1}
\right)
\\ &= t'  (\sigma  (\tau ^j)^{-1}  x) (e_0)= t'( (\tau ^j)^{-1} \sigma   x)
(e_0) = t'(\sigma    x) (e_0 \tau ^j).
\end{array}
$}

\noindent
Precise that   $\pi $ acts diagonally on $B( N^*X)^{\otimes p}
\otimes  B( N^*X)^{\otimes p}$  and on $W \otimes W$  and
that $\Delta $ (resp. $ \psi$) are $\pi $-linear. We
check in detail that:
$
 D't' (x) (w) = (t'\cup t')(x)(w)$, $ w \in W_0 $, $  x \in \overline{B(
N^*X)^{\otimes p} }
$.
Once again, it suffices to prove this formula  for the elements $ e_0
\sigma \in W_0$:

\centerline{$
\begin{array}{ll}
D't' (x)(e_0 \sigma  ) &= D t'(x) (e_0 \sigma ) + t'(dx)(e_0 \sigma )
= d(t'(x)(e_0 \sigma )) + t'(x)(\partial e_0 \sigma ) - t'(dx)(e_0 \sigma )\\
&= d(t(\sigma ^{-1} x))  + t(\sigma ^{-1} dx)
= (Dt) (\sigma ^{-1} x)
=(t \cup t) (\sigma ^{-1} x)\\
&= \left( m_{N^*X} \circ (t\otimes t) \circ \Delta  \right) (\sigma ^{-1} x)
=  m_{N^*X} \circ (t\otimes t)( \sigma ^{-1}  \Delta   x)
\\
&
=  m_{N^*X} \left(  t'(\sigma ^{-1} x')  (e_0) \otimes  t'(\sigma ^{-1}
x'')(e_0) \right)
=  m_{N^*X} \left(  t'(x')   (e_0\sigma ) \otimes  t'(x'')  (e_0 \sigma )
\right) \\
&= \left( m_{N^*X} \circ (t'\otimes t') \circ \Delta \right)   ( x) (\psi( e_0
\sigma ))
= (t'\cup t') (x) (e_0\sigma )\,.
 \end{array}
$}

Secondly,   for an arbitrary generator  $e_i \in W_i$, $i\geq 1$ and for
$[x] \in B_1( (N^*X)^{\otimes p} )$, we define
\quad
$
(\star _1) \qquad t'([x]) (e_i \sigma ) = \left\{
\begin{array}{lr}
0 & \mbox{if } x\  \in V_k \\
 \lambda (\sigma ^{-1} x)(e_i) &\mbox{if not }
\end{array}\,,
\right.
$
 where $V_k = {\Bbb F}_p
^{\otimes p-k-1 }\otimes N^*X \otimes   {\Bbb F}_p
^{\otimes p-k }$and $\lambda $ is a natural $\pi$-linear
transformation of degree 0 such that the following diagram commutes up to
homotopy

\centerline{$
\begin{array}{ccc}
   ( N^*X)^{\otimes p} &\stackrel{\lambda _X }{\rightarrow }&\mbox{Hom} (W,
N^*X)\\
   & m^{(p)}_{N^*X}\searrow & \downarrow ev_0\\
   & & N^*X.
\end{array}
$}

\noindent
Such a transformation  exists for  Corollary B.4 applies.  Observe  that if
$w \in W_i$ and $x \in  ((N^*X)^{\otimes p })$
then
$
\left (d_{N^*X} (\lambda _X (x)) - \lambda _X (dx) ^{\, }\right)(w)
=(-1)^{|x|}
\lambda(x)(\partial w)\,. $ The  formula ($\star _1$)  defines a $\pi
$-linear map
 $t' : {B_1\left( (N^*X)^{\otimes p} \right) } \to \mbox{Hom} (W, N^*X) $
and since  $
\overline \Delta [x] =0$ it
follows that
$
Dt' ([x]) =0= (t'\cup t')[x] \,, \qquad [x] \in {B_1\left((N^*X)^{\otimes
p}\right)}\,.
$

\vspace{2mm}

Thirdly, we assume that, for some $k\geq 1$, there is a $\pi $-linear map
$t' : \overline{B_{<k}\left( (N^*X)^{\otimes p} \right) } \to \mbox{Hom}
(W, N^*X) \,,
$
which extends to a $\pi $-linear map
$
t' : B_k\left( (N^*X)^{\otimes p} \right) \to \mbox{Hom} (W _{<i} , N^*X) \,,
$
for some $i \geq 1$ and  satisfying $ D't'(y) = (t'\cup t')(y)\,, \quad y
\in \overline{B_{<k}\left(
(N^*X)^{\otimes p} \right) }$. Let $d _1$  be the internal part of the
differential in the bar construction $d
= d_1+ d_2 $. The contravariant functor
$
X \mapsto FX =  ({\Bbb F}_p \oplus s^i B_k\left( (N^*X)^{\otimes p} \right)
, d)
$
where
$d(s^i y)= s^i d y \,, \quad y \in  B_k( (N^*X)^{\otimes p})$ is acyclic on
the
models $\{ \triangle ^n \}_{n\geq 0}$. We  define the natural transformation:
$
{\cal T}_{i, X}  : FX \to N^*X \,, \qquad y \mapsto  {\cal T}_{i, X} (y)=
(t'\cup t')(y)(e_{i}) - (-1) ^{|y|} t'(y) (
\partial e_{i} ) -t'(d ^2y)(e_{i})
\,.
$ (The inductive hypothesis ensures us that ${\cal T }_{i, X} (y)$ is well
defined.)
Let us check in detail  that $\tilde D  {\cal T}_{i, X} =0 $
when $\tilde D$ denotes the differential in $\mbox{Hom} (FX, N^*X)$.
Observe that, for   any
generator  $e_i $ of $W_i$ and if $y \in FX$,
the inductive hypothesis  implies :  $(D't')(y)(\partial e_i) =
(t'\cup t') (\partial e_i)$ and
$(D' t' ) (d_2y)(e_i) = (t'\cup t')(d_2y)(e_i)$. Thus  we obtain:

\centerline{$
\begin{array}{ll}
\left( \tilde D  {\cal T}_{i, X} \right) (y) &=
 d _{N^*X} ((t'\cup t')(y)(e_i) - (-1)^{|y|}  (D't')(y)(\partial e_i) \\
& \hspace{ 1 cm} + (-1) ^{|t'(y)|} t'(y) (\partial ^2 e_i) + (-1)^{|t'(y)|}
t'((d_1+d_2))(y))(\partial e_i)\\
& \hspace{ 1 cm} - ((D't')(d_2y)(e_i) + (-1)^{|t'(D_2(y)|} t'(d_2y)(\partial
e_i)
 + (-1)^{|t'|} t' ((d_1+d_2)(d_2)y) (e_i) ) \\
& \hspace{ 1 cm} -  \left( (t'\cup t')(d_1y)(e_i) -
(-1) ^{|y|+1} t'(d_1y)(\partial e_i) - t'(D_2D_1 y)(e_i) \right)\\
 &=  D'(t'\cup t')(y)(e_i) + (-1)^{|y|}( t'(d_1+d_2)(y)(\partial e_i) -
t'(d_1(y))(\partial
e_i) - t'(d_2y))(\partial e_i))\\
& \hspace{1 cm} + (t'(d_1\circ d_2(y))(e_i) + t'( d_2\circ d_1 ( y))(e_i) )\,.
\end{array}
$}

\noindent
Now  $\left( \tilde D  {\cal T}_{i, X} \right) (y)=0$, since each line in
the last  equality  is 0.

 Theorem B.4, establishes the existence of a
natutal transformation
${\cal S } _{i, X}  : FX \to N^*X $ such that $D'{\cal S } _{i, X}
 = {\cal T } _{i, X} $ that is:  if $z \in FX$ then
$
 {\cal T} _{i, X}  (z) = d {\cal S} _{i, X} (z) -
 (- 1) ^{|{\cal S } _{i, X}| } d ^1 z \,.
$ We set  for $ y \in B_k\left( (N^*X)^{\otimes p} \right) $,
$
t'(y)(e_i \sigma ) = {\cal S } _{i, X} (\sigma ^{-1} y) \,.
$ This formula defines  a $\pi $-linear map
$
t' : B_k\left( (N^*X)^{\otimes p} \right) \to \mbox{Hom} (W _{<i+1} , N^*X)
\,. $ It remains to check that  $ D't' = t' \cup
t'$.   Let $e_i $ be any generator of
$W_i$ and $y \in FX$,

\centerline{$
\begin{array}{lll}
(t'\cup t')(y)(e_i) &= {\cal T} _{i, X}(s^iy) + (-1)^{|y|} t'(y)(\partial
e_i) +
t'(d_2y)(e_i)\\
&= d{\cal S } _{i, X}(s^iy) - (_1) ^{{\cal S } _{i, X}|} {\cal S } _{i, X}(s^i
d_1 y) +
(-1)^{|y]} t'(y)(\partial e_i) + t'(d_2y)(e_i)\\
&= d(t'(y)(e_i)) + t'(d_1+d_2)(y))(e_i) + (-1) ^{|y|} t'(y)(\partial e_i)\\
&=( D' t')(y)(e_i)\,.
\end{array}
$}

We have thus constructed a twisting cochain
$
t'_X \in \mbox{Hom} ^1 \left(B \left((N^*X)^{\otimes p}\right),
\mbox{Hom}(W, N^*X)\right) $
such that
 $t'_X$ is natural in $X$ and  $t'_X$ is $\pi $-linear.

\vskip 3mm
\noindent {\bf 2.3} { \bf Proof of part 2.}  Recall the notation introduced in
\cite{[BT]}-Part II-$\S 3$. In particular,
consider  the simplicial model
$K$ of
$S ^1$,  and the cosimplicial space $\underline X$ defined as $ \underline
X(n)=
Map (K(n), X)$  whose geometric realization
$ || \underline X|| $ is homemorphic to $X ^{S ^1}$. Moreover, there is a
natural equivalence  $
\psi _X  : Tot  C^* \underline  X   {\longrightarrow }   N^*||\underline
X||\,, $
where $TotC^* \underline X$ denotes the total complex of the simplicial
complex
$C^* \underline X$.
On the other hand,  J. Jones, \cite{[Jo]}-6.3, has proved that  there is
natural chain map
$
\theta  _X  : {\frak C} N ^*X  \to  Tot
C^* \underline X
$
such that, if $X$ is  1-connected,  the composite
$J_X = \psi _X \circ \theta _X $ induces an isomorphism
$ HH_*N^*X \cong  H^*X ^{S ^1}\,.$  We obtain  the following diagram in the
category of
$\pi$-differential graded modules:

\centerline{$
\begin{array}{cccccc}
{\frak C} \left(\Omega B (N^*X)^{\otimes p}\right)  & \stackrel { \phi \circ
{\frak C}  \tilde
\kappa _X }
 {\longrightarrow} & \mbox{Hom} (W, {\frak C}N^*X)\\
 {\frak C}( \alpha _{(N^*X)^{\otimes p}}) \downarrow & & \hspace{1cm}
\downarrow
\mbox{Hom}(W, \theta _X \circ \psi_X)   & \\
{\frak C}\left( (N^*X)^{\otimes p}\right) &&\mbox{Hom}(W, N^*||\underline X
|| ) &\\
sh ^{(p)} \uparrow && \uparrow \tilde \Gamma _X\\
({\frak C} N^*X )^{\otimes p}  &\stackrel{\psi _X ^{\otimes p}}{
\longrightarrow}  &  (N^*||\underline X ||)^{\otimes p}
\end{array}
$}

\noindent
where $  \Gamma _X  $ denotes any natural structural map  defining
the usual Steenrood operations,  \cite{[Ma]}.  The functor $  X \mapsto
{\frak C}\left( (N^*X)^{\otimes p }\right)$
preserves the units and  is acyclic and  the
functor
$ X \mapsto  \mbox{Hom} (W,N^*||X||) $ is  corepresentable on  the models
$\cal M$, (see \cite{[BT]}-II-lemma 3.9) . Therefore, by Corollary  B.4,
there exists a $\pi
$-linear  natural transformation
$
 T_X : {\frak C}\left( (N^*X)^{\otimes p} \right)  \to  \mbox{Hom}
(W,N^*||X||)$
 such that
$T_X \circ {\frak C} \alpha _{(N^*X)^{\otimes p}} \simeq _{\pi}
\mbox{Hom}(W,\theta _X
\circ \psi_X)\circ  \phi \circ {\frak C}\tilde
\kappa _X
$ and $T_X \circ sh ^{(p)} \simeq _{\pi} \tilde \Gamma _X \circ \psi_X
^{\otimes
p}$. Consequently,
the second part  of theorem B is proved.

\vspace{1cm}

\centerline{ \bf \large $\S3$   $\pi$-$shc$  models.}

\vspace{3mm}

We continue the study of the homotopy category of
$\pi$-$shc$ algebras by constructing a convenient  model of these
algebras. As
it would be expected this model is rather
complicated.  Nevertheless, in the topological case, it provides  a
theoretical
object simpler than the singular cochain
algebra.

\vspace{3mm}
       \noindent {\bf 3.1 Minimal algebra.}  Let $V =\{V^i\} _{i \geq 1} $
be  a   graded vector
space and let  $TV$ denotes  the free graded
algebra generated by $V$:  $T^rV = V \otimes V \otimes ... \otimes V $
($r$-times)  and
$v_1v_2...v_k \in  (TV)^n $ if $\sum_{i=1}^k  |v_i| =n$.  The  differential
$d$ on
$TV$ is the  unique degree 1 derivation on
$TV$ defined by a given linear map $ V \to TV$ and such that $d_V \circ d_V
=0$.
 The
differential $d_V : TV \to TV$ decomposes as $d= d_0 +d_1+....$ with $d_kV
\subset T^{k+1}V$. In particular  $(V, d_0)$ is a
differential graded vector space. If we assume that $V^1=0$ and $d_0=0$
then  $(TV,d_V)$ is a called a
{\it 1-connected minimal algebra}. For any differential graded algebra
$(TU,d_U)$
such that $H^0(TU,d_U)= {\Bbb F}_p$  and $H^1(TU,d_U)=0$ there exists a
sequence
of homomorphisms  of differential graded algebras,

\centerline{$
(TU,d_U) \stackrel{p_V} \rightarrow (TV, d_V)  \stackrel{\varphi
_V}\rightarrow
(TU,d_U) $}

\noindent
where
$(TV,d_V)$ denotes  a  1-connected minimal algebra,  $p_V \circ \varphi _V
= id$,  $
\varphi _V \circ p_V  \simeq _{\bf DA}  id$ and $V$ such that $V\cong H
(U,d_{U,0})$. Moreover,   \cite{[HL]},
\cite{[FHT]}, \cite{[FHT2]}, $(TV,d_V)$  is unique up to  isomorphisms.

\vspace{3mm}
       \noindent {\bf 3.2 Minimal model of a product.}    Assume that $(A,
d_A)$ is a differential graded algebra such that
$H^0(A, d_A) = {\Bbb F}_p
$ and
$ H^1(A, d_A) = 0$, and let
$
 (T U[n],
 d_{ U[n] }) = \Omega \left( (BA)^{\otimes n} \right)$, $n\geq 1$.  By the
discussion above, for each $n \geq 1$,  we obtain
a sequences:

\centerline{$(T U[n],
 d_{ U[n] }) = \Omega \left( (BA)^{\otimes n} \right)\stackrel{ p_{ V[n]}}
\rightarrow (T V[n],
 d_{ V[n] })\stackrel {\varphi _{V[n]}}\rightarrow(T U[n],
 d_{ U[n] })$}

\noindent
 with
 $
V[n]  = s^{-1} \left( \overline{H \left( (BA)^{\otimes n} \right) } \right)
\cong
s^{-1} \left( \overline{ \left(  H  (BA) \right) ^{\otimes n }}\right)  =
s^{-1} \left( \overline{\left( { \Bbb F}_p \oplus s V \right)^{\otimes n}
}\right)  $\\
 $ \cong  \left( \bigoplus _{k=1} ^n \left( \left( { \Bbb F}_p
\right)^{\otimes k-1 } \otimes V \otimes  \left(  { \Bbb F}_p
\right)^{\otimes n-k } \right) \right)
\oplus ... \oplus s^{-1} \left( sV \otimes sV\otimes ...\otimes sV \right) $.

For $n=1$,  $V[1] = V =  s^{-1}\overline{ H(BA)} $
 and   the  composite
$
 \psi _V =\alpha _A \circ \varphi _V :  (TV,d_V) \to  A
$
is a
quasi-isomorphism. The algebra   $(TV, d_V)$
is called the {\it 1-connected minimal model} of $A$.

For $n\geq 2$, consider the homomorphism
$q_{V[n]}: (T V[n], d_{ V[n] })  \to  (TV, d_V)^{\otimes n}$ defined by
$ q_{V[n]}(y)=  1 ^{\otimes k-1} \otimes y \otimes 1^{\otimes n-k}$ if $y
\in  V_k := { \Bbb F}_p
^{\otimes k-1 }\otimes V \otimes   { \Bbb F}_p
^{\otimes n-k }$, $k=1,2,...,n$ et $q_{V[n]}(y)= 0 $ if $ y \in V[n]-
\bigoplus_{i=1}^n  V_i$.  The
composite
 $(T V[n],
 d_{ V[n] })  \stackrel{q_{V[n]}}\rightarrow (TV, d_V)^{\otimes
n}\stackrel{(\psi _V)^{\otimes n}}
\rightarrow
(A, d_A)^{\otimes n})$ is a  quasi-isomorphism  (\cite {[BT]}-6.5).
 Therefore $(T V[n],
 d_{ V[n] })$ is the  minimal model of $A^{\otimes n}$.

\vspace{ 3mm}
\noindent{ \bf 3.3 Lemma.} (Equivariant lifting lemma.) {\it Let $\pi$ any
finite group and consider

i)   a  graded  free $\pi $-module, $U = \{U^i \}_{i\geq 2}$,

ii)   a   minimal model,  $(TU,d)$,  with $\pi $ acting by the rule $
\sigma \cdot
u_1u_2...u_k = \sigma u_1...\sigma u_k$, so that $(TU,d)$ is a $\pi $-complex,

iii) two differential graded algebras,  $A$ and $B$,  which are also $\pi
$-complexes,

iv)   a homomorphism of differential graded algebras, $f : (TU,d) \to B $,
which is $\pi
$-equivariant,

v) a sujective quasi-isomorphism,  $\varphi : A \to B$,  of differential
graded algebras  which is
$\pi $-equivariant.

\noindent Then there exists a homomorphism of differential graded algebras
$ g : (TU,d) \to B $
which is $\pi
$-equivariant  and such that $\varphi \circ g = f$.}

\vspace{ 3mm}
\noindent{ \bf Proof.}
 We choose a homogeneous linear basis $\{ u_i \}_{i \in
 I} $  of the $\pi$-free graded module  $U$ with $I$ a well-ordered set.   We
denote by  $ U_{<i}$ the sub-$\pi $-module generated by the elements $ u_i$,
$j<i$. The minimality
condition together with the 1-connectivity condition imply that
$ d u_i \in T(U_{<i})$. Suppose that $g : ( T(U_{<i}), d) \to B$ has been
constructed so that $g$ is
$\pi
$-equivariant  and  $\varphi \circ g = f$. Then $g(du_i) \in B \cap \mbox{
Ker } d_B$ and
$\varphi \circ g (du_i) = d_A f(e_i)$. Therefore there exists $ b_i \in B $
such
that $g(du_i)=
d_Bb_i$. Moreover, $d_Af(u_i) = d _A \varphi b_i$. For $\sigma \in \pi$, we
set
$g(\sigma u_i )=
\sigma (b_i+z_i)$, and we check that $\varphi \circ g(\sigma u_i)= f(\sigma
u_i)$ and $d
g(\sigma u_i) = g(d\sigma u_i) $.

\vskip 4mm
\noindent {\bf 3.4 $\pi$-$shc$-minimal model. } For any integer $n>1$ the
group ${\frak S}_n$  acts on
$V[n]\subset s^{-1} \left( H  (BA) \right) ^{\otimes n}$.  This action extends
diagonally on $TV[n] $ so that the differential $d_{V[n]} $ and the
homomorphism
$(\psi _V)^{\otimes n}\circ q_{V[n]}$
 are ${\frak S}_n$-linear.  Since the natural  map $ \alpha_{A^{\otimes n}}
$ is  a  ${\frak S}_n$-equivariant
surjective
quasi-isomorphism,
 Lemma 3.3 implies that the composite $ (\psi _V)^{\otimes p}\circ
q_{V[n]}$ lifts
 to a
homomorphism of differential graded algebras
$
 L : TV[n]
  \to \Omega    B (A^{\otimes n})
$
which is ${\frak S}_n$-equivariant and
$ \alpha_{A^{\otimes n}}\circ L = (\psi _V)^{\otimes n}\circ q_{V[n]}$.
Let $(A,d_A, \mu_A)$ be an augmented
$shc$ algebra and
assume that
$H^0(A, d_A) = {\Bbb F}_p $ and that $
       H^1(A, d_A) = 0$.  The   composite, $ \mu _V ^{(n)} := p_V \circ \mu
_A ^(n) \circ L :
(TV[n], d_{V[n]}) \to   (T(V),d_V)
$
 is a homomorphism of differential graded algebras.
  If $n=2$, $ \mu _V ^{(2)}: = \mu
_V :  (T( V[2]) ,  d_{ V[2]}) \to (TV,d)$. The triple
$(TV, d_V, \mu _V)$ is called a {\it $shc$-minimal model for $(A, d_A, \mu
_A)$} (See \cite{[BT]}-I-$\S6$ for
more details.)

    Let  $(A, d_A, \mu_A, \tilde \kappa_A)$ be a
$\pi$-$shc$  algebra such that $H^0A= {\Bbb F}_p$ and $H^1A=0$. By Lemma
3.3,  $\tilde
\kappa_A \circ L $ lifts to a $\pi $-equivariant homomorphism of algebras
 $\hat   \kappa_A : TV[p] \to \mbox{Hom}(W, \Omega   B A)$.
Since   $\pi $ is assumed to act trivially on $TV $  the composite
 $
\tilde \kappa_V = \mbox{Hom}(W, p_V)\circ \hat \kappa_A :
TV[p] \rightarrow \mbox{Hom}(W, TV)
$
 is a $\pi$-equivariant map and we obtain the commutative   diagram:

\centerline{$
 \begin{array}{cccc}
& \Omega B(A^{ \otimes p})  & \stackrel{\tilde \kappa_A}{\longrightarrow} &
\mbox{Hom}(W, A)
\\
&\Big\uparrow {L}   & &\hspace{20mm}\Big\uparrow \mbox{Hom}(W, \alpha _A)
\\
& (TV[p], d_{V[p]})  & \stackrel{\hat \kappa_A}{\longrightarrow} &
\mbox{Hom}(W, \Omega BA)
\\
&  & \hspace{-15mm} \tilde \kappa _V \searrow  &\hspace{20mm}\Big\downarrow
\mbox{Hom}(W, p _V)
\\
&  & & \mbox{Hom}(W, TV)\,.
\\
\end{array}
$}

\noindent
The triple  $(TV, d_V, \mu _V, \tilde \kappa _V )$ is called a {\it
$\pi$-$shc$  minimal
model for the $\pi$-$shc$ algebra  $(A, d_A,  \mu _A, \tilde \kappa _A )$}.

\vspace{3mm}
       \noindent {\bf 3.5. Proposition.} {\it
 Let $(A,d_A, \mu_A,\tilde \kappa _A )$  and  $(TV, d_V, \mu _V, \tilde
\kappa _V )$ be
      as above  then:

1) $(TV, d_V, \mu _V, \tilde \kappa _V )$ is a $\pi$-$shc$ algebra,

2) $\psi _V : (TV,d) \to A $ is a strict $\pi$-$shc$ homomorphism and a
quasi-isomorphism,

3) the canonical maps  $
(TV,d_V) \stackrel {\iota _V } {\longrightarrow} {\frak C}(TV,d_V) \stackrel
{\rho _V } {\longrightarrow} B(TV,d_V)
$  respect the Steenrood operation in homology.

4)  the quasi-isomorphism $\psi _V : (TV,d_V) \to (A,d_A) $ induces a
commutative diagram

\centerline{$
\begin{array}{cccccc}
(TV,d_V) &\stackrel {\iota _V } {\longrightarrow} {\frak C}&(TV,d_V) \stackrel
{\rho _V } {\longrightarrow} &B(TV,d_V) \\
\psi_V \downarrow & & \downarrow {\frak C}\psi_V&   \downarrow  B\psi_V \\
(A,d_A) &\stackrel {\iota _A } {\longrightarrow} {\frak C}&(A,d_A) \stackrel
{\rho _A } {\longrightarrow} & B(A,d_A)
\end{array}
$}

\noindent
in which all vertical arrows
induce isomomorphisms in homomology which respect the Steenrod operations.
}

\vspace{2mm}
\noindent{\bf Proof.}  From the definition  of $\tilde \kappa_V$,  we deduce
that
$ev_0\circ \hat \kappa_V\simeq_{{\bf DA}}\mu^n_V$.
This implies that
$(TV, \mu_V,  \tilde \kappa_V)$ is a $\pi$-$shc$ algebra.  Since
$p_V\circ\varphi_V=id_{TV}$,
$\varphi_V\circ p_V
\simeq _{\bf DA} id _{\Omega BA}$    then $Hom(W, \phi _V ) \circ \tilde
\kappa _V \simeq _{\pi \! -\!{\bf DA}} \hat \kappa
_A$. Indeed, it is straighforward to check that  if
$\pi$ acts trivially on $A$ and $A'$ and if
  $f, g: A\rightarrow A'$ are homotopic  in {\bf DA},
  then $\mbox{Hom}(W, f)$ and $\mbox{Hom}(W, g)$ are $\pi$-homotopic in
{\bf DA}. It remains to show
 that $\psi _V$ is a homomorphism $\pi$-$shc$ algebras. This follows
directly from the commutativity up to $\pi$-
homotopy in {\bf DA} of the diagram

\centerline{$
\begin{array}{ccccccccccc}
 \Omega \left( (B A)^{\otimes p} \right) & \stackrel{\tilde \kappa_A }
 {\longrightarrow} & \mbox{Hom} (W, A) \\
\hspace{-7 mm} L \uparrow && \hspace{15mm} \uparrow \mbox{Hom} (W, \psi _V)\\
T  V[p] & \stackrel{\tilde \kappa _V }{\longrightarrow} & \mbox{Hom} (W, TV)\,.
\end{array}
$}

\vspace{3mm}
       \noindent {\bf 3.6. Proposition.} {\it  Let $(A,d_A, \mu_A,\tilde
\kappa _A )$  and  $(TV, d_V, \mu _V, \tilde \kappa _V )$ be
      as above.   Assume $p=2$ and   set $\hat V = V[2]$, \quad $V' = V
\otimes {\Bbb
F}_p$,\quad $V'' = {\Bbb F}_p\otimes V$ \quad and  $  V' \# V'' =
s^{-1} (sV \otimes sV) $,  so that:
$\hat V = V'\oplus V''\oplus V' \# V''\,.$. If  $\tilde  \kappa:\Omega
B(A\otimes A) \to \mbox{\rm Hom}(W, A)$ satisfies:

\centerline{$
(\star) \quad \tilde  \kappa (\langle [a\otimes b]\rangle ) = 0\mbox{ if }
[a\otimes b] \in B_1(A\otimes A)
  \mbox{ and } \epsilon_A(a)=0\mbox{ or } \epsilon_A(b)=0 \,,
$}

\noindent
 then,  for any $ v\in
V^n\cap kerd_V$,
$
Sq^i(cl(v)) = \tilde  \kappa _V(v\#v)(e_{n-i-1}) \,.
$}

\vspace{2mm}

 Observe that if $A = N^*X$, condition $(\star )$  is  satisfied
(see 2.2, ($\star_1$)).

\vspace{2mm}
\noindent{\bf Proof.} Define
$
 K_A : \Omega (BA\otimes BA)  \rightarrow \Omega BA\otimes \Omega BA
$
 as follows: we write $\Omega (BA\otimes BA) = (T(s^{-1}(B^+A\otimes {\Bbb
F}_p)\oplus s^{-1}({\Bbb F}_p\otimes B^+A)\oplus
s^{-1}(B^+A\otimes B^+A)), D) $ and we  set $
K_A (s^{-1}[a_1\vert ...\vert a_i]\otimes 1) = s^{-1}[a_1\vert ...\vert
a_i]\otimes 1 $, $
 K_A(1\otimes s^{-1}[a_1\vert ...\vert a_i]) =  1\otimes
s^{-1}[a_1\vert...\vert a_i]$, $ K_A(s^{-1}([a_1\vert ...\vert
a_i]\otimes [b_1\vert ...\vert b_j]) ) = 0 \,, i,j >0$.
 One   easily check (\cite{[BT]}-I-6.4 )  that
$K_A$ commutes  with the
differentials. Since, $(\alpha _A\otimes \alpha _A)\circ K_A = \alpha
_{A\otimes A} \circ \Omega \overline{sh}$, the
homomorphism $K_A$ is a
surjective quasi-isomorphism. The map
$
 S: W\otimes_{\pi} (A\otimes A)  \longrightarrow
W\otimes_{\pi}\Omega(BA\otimes BA)
$
  defined  by  $
       S(e_0\otimes a\otimes b ) = e_0\otimes \langle [a]\otimes 1\rangle
\langle 1\otimes [b]\rangle $ and if  $i>0$ by  $
  S(e_i\otimes a\otimes b ) = e_i\otimes\langle [a]\otimes 1\rangle \,
       \langle 1\otimes [b]\rangle    + e_{i-1}\otimes\langle [b]\otimes
        [a]\rangle $
satisfies:  $id_{W} \otimes  (\alpha_A\otimes \alpha_A )
  \circ id_{W}\otimes K_A  \circ S =
  id_{W \otimes  A\otimes A}$ and  $Sd = dS  $.
  As in 1.6,  we deduce:

\centerline{$
\begin{array}{rl}
Sq^ix& = cl(\tilde  \kappa (i_{A^{\otimes 2}}(a\otimes a))(e_{n-i}))
=  cl( \kappa (e_{n-i}\otimes i_{A^{\otimes 2}}(a\otimes a)))
= cl(\kappa_{/\pi }(e_{n-i}\otimes i_{A^{\otimes 2}}(a\otimes a))) \\
&=  cl(\kappa_{/\pi }(((id\otimes \Omega \overline {sh} )
\circ S)(e_{n-i}\otimes i_{A^{\otimes 2}}(a\otimes a))))
\\
&= cl(\kappa_{/\pi }(e_{n-i}\otimes \langle [a\otimes 1]\rangle.
       \langle [1\otimes a]\rangle ))   + cl(\kappa_{/\pi }( e_{n-i-1}\otimes
       \Omega \overline {sh} (\langle [a]\otimes
        [a]\rangle )))
\\
&=cl(\tilde \kappa(  \Omega \overline {sh} (\langle [a]\otimes
        [a]\rangle)) ( e_{n-i-1}))
\,.
\end{array}
$}

\noindent
The last equality comes from the condition $(\star) $ and the fact that $
\Omega \overline {sh}$ is
$\pi$-linear.

Now observe that $ \varphi_{ \hat V}: T(  V'\oplus V''\oplus V' \# V'')
\rightarrow T(s^{-1}
\overline {BA \otimes BA})$ identifies $v'\in V'$ with $[a]\otimes 1$,\quad
$v''\in V''$ with $1\otimes [a]$
and $v'\#v''\in V' \# V''$ with $\langle [a]\otimes
        [a]\rangle$. Thus if $v\in V^n\cap kerd_V$ and $x=cl(v)$ then
$H\psi_V\circ Sq^i (cl(v)) =
 cl(\tilde \kappa (\Omega \overline {sh} \circ \varphi_{ \hat
V}(v'\#v''))).$ The formula now  follows from the definition of
$\tilde
\kappa_V$.

  \vspace{3mm}
       \noindent {\bf 3.7 $shc$-equivalence.} Two $shc$ algebras (resp.
$\pi$-$shc$ algebras) $A$ and
$A'$  are {\it $shc$
       equivalent}, $A \simeq_{shc} A'$ (resp. {\it $\pi$-$shc$
       equivalent}, $A \simeq_{\pi\mbox{-}shc} A'$) if there exists a
sequence of strict $shc$ (resp.  $\pi$-$shc$)
homomorphism  $A \leftarrow A_1\rightarrow ......\rightarrow A'$ inducing
isomorphisms in homology.  If $A
\simeq_{\pi\mbox{-}shc} A'$, then  $H(A) \cong  H(A')$, and $A$, $A'$ have
the same  $\pi$-$shc$ minimal model. Two spaces
$X$ and $Y$ are $shc$ equivalent (resp. $\pi$-$shc$ equivalent) if the
differential graded algebras $N^*X$ and
$N^*Y$ are $shc$ equivalent (resp.
$\pi$-$shc$ equivalent).

\vspace{3mm}
       \noindent {\bf 3.8. Example.} We  exhibit  two spaces $X$ and $Y$ which
are  $shc$ equivalent but not   $\pi$-$shc$  equivalent.    Let  us consider
$X= \Sigma ^2 {\Bbb C}P ^2$ and $Y = S^4\vee S^6$. The  spaces $X$ and $Y$
have the same
        Adams-Hilton model  namely $(T(x_3, x_5), 0)$. This shows, in
particular, that
        $H_*\Omega X \cong H_*\Omega Y $, as graded Hopf algebras, and that
$X$
and $Y$ are ${\Bbb F}_2$-formal. Furthermore,
        $H^* =H^*X = {\Bbb F}_2a_4 \oplus {\Bbb F}_2a_6 \cong H^*Y$ as graded
algebras (with trivial products). Hence $X$ and $Y$ have the same minimal
model say $(TV,d_V)$ with
        $V = s^{-1}H_+\Omega X = s^{-1}T^+(x_3, x_5)$,
        $ d_Vs^{-1}x_3=  d_Vs^{-1}x_5=0$. The map $\psi: (TV,
d_V)\rightarrow H^*$ defined by
$\psi (s^{-1}x_3) = a_4,\quad
       \psi (s^{-1}x_5) = a_6  \mbox{ and } \psi = 0  \mbox{ on }
       s^{-1}T^{\geq 2}(x_3, x_5)$
       is a surjective quasi-isomorphism.  Let $(TV, d_V, \mu_1, \kappa
_1)$ (resp. $(TV, d_V, \mu_2, \kappa _2)$ be a
$\pi$-$shc$ minimal model of $X$ (resp. $Y$), then for $i=1, 2$
       $\mu_i: (T{\hat V}, d_{\hat V})\rightarrow (TV, d_V)$
       identifies
       the two copies of $V$. Moreover,
       $\vert \mu_i(x\#y)\vert >6$ for any  $x, y \in V$.
        Therefore
       $\psi\circ \mu_i(x\#y) = 0$.
This yields the commutative diagram

\centerline{$
\begin{array}{ccccc}
T{\hat V}&
  \stackrel{\psi\otimes \psi\circ q_{\hat V} }{\longrightarrow} &
 H^*\otimes H^*&
 \stackrel{\psi\otimes \psi\circ q_{\hat V} }{\longleftarrow}&
  T{\hat V}\\
   \mu_1\Big\downarrow &
&\Big\downarrow m_{H^*}& &\Big\downarrow \mu_2\\
 TV &\stackrel{\psi}{\longrightarrow} &
  H^*&\stackrel{\psi}{\longleftarrow} &TV\\
\end{array}
$}

\noindent
 This proves that $X$ and
$Y$ are $shc$ equivalent. Now we prove that $X$ and $Y$ are not $\pi$-$shc$
equivalent. Recall that if $a_4 \in H^4  \Sigma
^2 {\Bbb C}P ^2$ then  $Sq^2a_4 = a_6$
 while  if $ a_4 \in H^4 (S^4\vee S^6)$ then  $Sq^2a_4 = 0$. Thus,  by
Proposition 3.6,
  $ \tilde \kappa_1(a_4\#a_4)(e_1)=a_6$ while   $ \tilde
\kappa_2(a_4\#a_4)(e_1)=0$.
  Now  $\tilde \kappa_1$,   $ \tilde\kappa_2$ are not $\pi$-homotopic, and
hence $X$ and $Y$ are not
$\pi$-$shc$ equivalent.

 \vspace{3mm}
       \noindent {\bf 3.9. Example.} We compute the Steenrod operations on
$HH_*A$ when $A = {\Bbb F}_2[u] = T(u)$,
$\vert u\vert =2$.

First, recall that, endowed with the shuffle product, $
 {\frak  C}T(u)  \stackrel{sh} \rightarrow  {\frak  C}(T(u) \otimes T(u))
\stackrel{{\frak  C}m} \rightarrow  {\frak
C}T(u)$, the complex
${\frak  C}T(u)$ is a commutative differential graded algebra. We consider
the homomorphism of differential graded algebras
$\rho: (T(u)\otimes \Lambda su, 0) \rightarrow {\frak  C}T(u)$
defined by
 $\rho (u\otimes 1) = u[\, ]$,  $\rho (1\otimes su) = 1[u]$. As proved in
\cite{[L]}-Proposition 3.1.2, $\rho$ is a
quasi-isomorphism of differential graded algebras. In particular,
$T(u)\otimes \Lambda su\cong HH_*T(u)$ as commutative
graded algebras.

Secondly, we define a structure of $\pi$-$shc$ algebra on $T(u)$. Let
$V=u{\Bbb F}_2$ and set ${\hat V}=u'{\Bbb
F}_2\oplus  u''{\Bbb F}_2\oplus u'\# u''{\Bbb F}_2$. We define $\tilde
\kappa: T({\hat V})\rightarrow \mbox{Hom}(W, T(u))$ by
$
\tilde \kappa  (u')(e_i)= \tilde \kappa  (u')(e_i\tau)= \tilde \kappa
(u'')(e_i)=\tilde\kappa  (u'')(e_i\tau )= u $ if $i=0$
and $=0$ if $i>0$ and
$
\tilde \kappa  (u'\#u'')(e_i) = \tilde \kappa  (u'\#u'')(e_i\tau )=
u $ if $  i=1 $, $= 0$ if $i>0$.
As   in the proof of Proposition 3.6, $\tilde \kappa $ defines the Steenrod
operations on $HT(u) = T(u)$:
 $
Sq^0(u)=
u$, $Sq^2 = u^2 $ and $Sq^i u =0 $ if $ i\neq 0, 2$. By Proposition 3.6 one
deduces that $\tilde \kappa$ is the unique
$\pi$-$shc$  structural map on $T(u)$ with the Steenrod operations defined
above.

Now in order to compute  the Steenrod operations on $HH_*T(u)$ we consider
the structural map
$ \tilde \theta= \phi \circ {\frak C}\tilde \kappa  : {\frak C} \left(
T{\hat V} \right)
\to \mbox{Hom}( W, {\frak C}T(u))$,    and a linear section $S$ of
$id_W\otimes{\frak C} q_{\hat V}:  W\otimes_{\pi}\frak{C}(T{\hat V}
)\rightarrow  W\otimes_{\pi}\frak{C}(T(u)\otimes T(u) )$.
Then Steenrod operations are defined by:
$
Sq^i (x) = cl( \theta (S\circ (id_W\otimes sh) (e_{n-i}\otimes x\otimes x)),
\quad x \in HH^nT(u)
$.

Finally observe that $S$ is uniquely determined in low degrees as follows:

\centerline{$
\begin{array}{lllll}
 S(e_i\otimes 1[\,])& =e_i\otimes 1[\,]   & S(e_i\otimes (u\otimes
1)[\,])&=e_i\otimes u'[\,]\\
S(e_i\otimes (1\otimes u)[\,]&=e_i\otimes u''[\,]& S(e_i\otimes 1[u\otimes
1])&=e_i\otimes 1[u'] \\
  S(e_i\otimes (u\otimes u)[\,])&=e_i\otimes u'u''[\,]+e_{i-1}\otimes
u'\#u''[ ] & S(e_i\otimes 1[1\otimes u])&=e_i\otimes 1[u'']\\
S(e_i\otimes (u\otimes 1)[u\otimes 1])&=e_i\otimes u'[u'] & S(e_i\otimes
(1\otimes u)[1\otimes u])&=e_i\otimes
u''[u'']\\
S(e_i\otimes (1\otimes u)[u\otimes 1])&=e_i\otimes u''[u'] + e_i\otimes
u'\#u''[\, ] &
\\ S(e_i\otimes (u\otimes 1)[1\otimes
u])&=e_i\otimes u'[u''] + e_i\otimes  u'\#u''[\, ] & \\
S(e_i\otimes 1[u\otimes u])&=e_i\otimes 1[u'u'']+e_{i-1}\otimes 1[u'\#u'']  &
S(e_i\otimes 1[u\otimes 1\vert 1\otimes  u])&=e_i\otimes 1[u'|u'']\\
S(e_i\otimes 1[1\otimes u\vert u\otimes  1])&=e_i\otimes
1[u''|u']+e_{i}\otimes 1[u'\#u'']
\end{array}
$}

\noindent
(Here we make the convention that $e_{-1}=0$.) Therefore, since  $\tilde
\theta$
respects the products and with the aid of
Proposition 3.6, we can do  the following computations:

\centerline{$
\begin{array}{ll}
Sq^0(u[\, ])& = cl(\theta\circ S\circ sh(e_2\otimes u[\, ]\otimes u[\, ]))
  =cl(\theta (e_2\otimes u'u''[\,] + e_1 \otimes  u' \#u''[\,]))
\\
&= cl(\tilde \theta ( u'u'')(e_2)[\,])) + cl(\tilde \theta ( (u' \#u'')
(e_1)[\,] ))
= cl( \tilde \kappa (u'\#u'')(e_1)[\,])
=
Sq^0 (u[\,]) = u[\,]
\\
Sq^1(u[\, ])& = cl(\theta\circ S\circ sh(e_1\otimes u[\, ]\otimes u[\, ]))
 =cl(\theta (e_1\otimes u'u''[\,] + e_0 \otimes 1[u' \#u'']))
\\
&= cl( \tilde \kappa (u'\#u'')(e_0)[\,]))
=0
\\
Sq^2(u[\, ])& = cl(\theta\circ S\circ sh(e_0\otimes u[\, ]\otimes u[\, ]))
= cl( \tilde \kappa (u'u'')(e_0)[\,]))
=u^2.
\\
Sq^0(1[u])& = cl(\theta\circ S\circ sh(e_1\otimes 1[u]\otimes 1[u]))
  =cl(\theta (e_1\otimes 1[u'\vert u''] + e_1 \otimes 1[u'\#u''] ))
\\
&= cl(1[\tilde \kappa (u'\#u'')(e_1)) = 1[Sq^0 u] =1[u]
\\
Sq^1(1[u])& = cl(\theta\circ S\circ sh(e_0\otimes 1[u]\otimes 1[u]))
  \\
  &=cl(\theta (e_0\otimes 1[u'\vert u''] + e_0\otimes 1[u''\vert u']+ e_0
\otimes 1[u'\#u''] ))
\\
&= cl(1[u\vert u] + 1[u\vert u]) + cl(1[\tilde \kappa (u'\#u'')(e_1)])
=
0.
\end{array}
$}

\noindent
Thus the Steenrod operations are completely defined on $HH_*T(u) =
T(u)\otimes \Lambda su$, by  the Cartan formula and the
formulas: $Sq^0  (u[\,]) = u[\, ]$, $Sq^2  (u[\,]) =u^2[\,]$, $Sq^i
(u[\,])=0$ for $i \neq 1,2$, $Sq^0  (1[u])=1[u]$,
$Sq^i  (1[u])=0$ if $ i>0 $. We recover  the topological Steenrod squares
 on $H^*\left(\left( {{\Bbb C}P ^{\infty}}\right) ^{S^1}\right) =H^*({\Bbb
C}P ^{\infty})
\otimes H^*(S^1)$. The  technics developped in this example,  can be
performed in order to study
the case when $A =  {\Bbb F}_2[u]/(u ^k)$,
$k\geq 2$.


\vspace{1cm }

\centerline{ \bf \large Appendix A -  Technicalities on $shc$-algebras.}

\vspace{3mm}
In  this appendix  we lay the material necessary  in order to complete the
proof of  Theorem A.
\vspace{3mm }

\noindent{\bf A.1 Lemma.} {\it Let $A$, $A'$ be two $shc$ algebras and let
$f$, $g$ $\in {\bf DA}(A,A')$ be such that
 $f \simeq_{{\bf DA}} g$. If  $f$ is a  $shc$ homomorphism  then $g$ is a
$shc$ homomorphism.}

\vspace{2mm }

\noindent{\bf Proof}  Let $\theta$ be a $(f, g)$-derivation such that
$f - g = d_{A'}\circ\theta +\theta\circ d_{A}$. Then $\theta$ defines a
coderivation $\theta': BA\rightarrow BA'$:
$\theta'([a_1\vert ....\vert a_k]) =
\sum_{i=1}^k(-1)^{\epsilon_i}
[fa_1\vert ..\vert fa_{i-1}\vert \theta a_i\vert...
\vert ga_{i+1}..\vert ga_k]
$
where $\epsilon_i=\sum_{j=1}^{i-1} |a_j| + (i-1)$
satisfying
$Bf - Bg = d_{BA'}\circ\theta ' +\theta ' \circ d_{BA}$.  Now let $f_1$,
$f_2 \in {\bf DC}(C, C')$ and consider  the
adjoint homomorphism  $ \tilde f_i   : C\to  C' \to B\Omega C'$  ($i=1,
2$).  If we assume that $f_1 \simeq_{{\bf DC}} f_2$
then
$\tilde f_1 \simeq_{{\bf  DC}} \tilde f_2$  and thus,  by \cite{[Mu]}-1.11,
$\Omega f_1 \simeq_{{\bf DA}} \Omega f_2$.
Therefore  if
$f \simeq_{{\bf DA}} g$ then $\Omega Bf \simeq_{{\bf DA}} \Omega Bg$. The
lemma follows now, from the obvious relations:
$
 \mu_{A'}\circ \Omega B(g\otimes g)\simeq_{{\bf DA}}
 \mu_{A'}\circ \Omega B(f\otimes f) \simeq_{{\bf DA}} \Omega Bf\circ \mu_A
\simeq_{{\bf DA}} \Omega Bg\circ \mu_A \,.
$

\vspace{3mm }

\noindent{\bf A.2 Trivialized  extensions.}  Let ${\bf TEX_A}$ be the
category of trivialized
extensions  of $A$ in {\bf DA}
in the sense of \cite{[Mu]}-2.1: $ X \stackrel{\alpha}{\to}A \in {\bf TEX}
$ if
$\alpha \in {\bf DA}(X,A)$ and if there exist $\rho \in {\bf DM}( A,X) $, $
h \in \mbox{Hom}(X,X)$ such that,
$\alpha \circ \rho = id_A$, $  \rho \circ \alpha -id_X = d_X \circ h + h
\circ d_X$, $
\rho \circ \eta _A = \eta _X$, $  \epsilon _X \circ \rho = \epsilon _A$, $
\alpha \circ h = 0$, $ h\circ \rho = 0$, $  h \circ h = 0$. A {\it (strict)
morphism of trivialized extension}, $f : ( X
\stackrel{\alpha}{\to}  A , \rho, h)
\to ( X' \stackrel{\alpha}{\to}  A , \rho', h')$  is a homomorphism $ f \in
{\bf DA}(X, X') $ such that $f \circ \rho = \rho
' $ and $f\circ h =  h'\circ f$.

 The following  facts are  proved in \cite{[Mu]} or  straightforward to prove.

 \vspace{2mm}
 \noindent{\bf F1} { If  $ X \stackrel{f}{\to} A \in {\bf TEX}_A$ and  $ Y
\stackrel{g}{\to} X \in {\bf TEX}_X$, then
$ Y \stackrel{f\circ g }{\to} A \in {\bf TEX}_A$ .}

 \vspace{2mm}
 \noindent{\bf F2}  { If  $ X \stackrel{f}{\to} A \in {\bf TEX}_A$ and  $
X' \stackrel{f'}{\to} A' \in {\bf TEX}_{A'}$,
then
$ X \otimes X'  \stackrel{f\otimes f'  }{\to} A \otimes A' \in {\bf
TEX}_{A\otimes A'} $ .}

 \vspace{2mm}
 \noindent{\bf F3}
 { If $A$ is a differential graded algebra  then $\Omega B A \stackrel
{\alpha _A}{\to} A$ in an initial object in $
{\bf TEX}$.}

 \vspace{2mm}
\noindent \noindent{\bf F4}  { If  $(A, \mu _A)$ and $(A', \mu_{A'})$ are
$shc$ algebras
then there exists a natural homomorphism $\mu_{A\otimes A'} $ such that $(
A\otimes A', \mu_{A\otimes A'}) $ is a $shc$
algebra. In particular, $A ^{\otimes n} \,, n \geq 2 $ is naturally a $shc$
algebra.}

\vspace{2mm}
\noindent{\bf F5} { The homomorphism  $\alpha_{\Omega BA}: \Omega B \Omega
BA\to \Omega BA$
admits  a section
$\theta_{\Omega BA}\in {\bf DA}(\Omega BA, \Omega B\Omega BA)$.}

\vspace{2mm}
\noindent{\bf F6} {  If $(A, \mu_A)$  is a $shc$ algebra  so is $\Omega B A
$.}

\vspace{2mm}
\noindent{\bf F7} {  $\Omega B\alpha_{A}\simeq_{{\bf DA}}\alpha_{\Omega BA}$.}

\vspace{2mm}

\noindent{\bf F8} { If $(A, \mu_A)$ is a $shc$ algebra  then $\alpha_{A} :
\Omega B A \to A $ is a strict
 $shc$  homomorphism.}

\vspace{3mm }

\noindent{\bf A.3 Lemma.} {\it  Let $(A,d,\mu )$ be a $shc$  algebra.
For any  $n \geq 2$,  there exists a homomorphism of
differential graded algebras (called the $shc$ iterated structural map)
$
\mu ^{(n)} : \Omega B \left( A ^{\otimes n} \right) \to \Omega BA
$
such that:  $ \mu ^{(2)} = \mu$ and  $\alpha _A \circ \mu ^{(n)} \circ i_{A
^{\otimes n}}
\simeq m_A  ^{(n)} $. Moreover,   $\mu^{(n)}$ is a strict homomorphism
of $shc$  algebras.
}

\vspace{2mm }

\noindent{\bf  Proof.}  Let $A_1$ and $A_2$ be two differential graded
algebras.  From F3, we deduce that there exists  a
natural homomorphism of differential graded algebras:
$
\theta_{_{A_1,A_2}} : \Omega B (A_1 \otimes A_2) \to
 \Omega B( A_1 \otimes \Omega B A_2)
$
such that:
$
(\star ) \quad  (id_{A_1} \otimes \alpha _{A_2} ) \circ
\alpha _{A_1 \otimes \Omega B A_2} \circ  \theta _{_{A_1, A_2}} = \alpha
_{A_1
\otimes A_2}\,.
$ It follows, from the unicity property , that if $A_3$ is a differential
graded algebra then $
(\theta _{_{A_1,A_2}} \otimes id_{\Omega B A_3 }) \circ \theta _{_{A_1 \otimes
A_2, A_3}}
= (id_{\Omega B A_1 } \otimes \theta _{_{A_2,A_3}}) \circ \theta_{ _{A_1,
A_2\otimes A_3}}\,.$  Suppose inductively that  $\mu ^{(n)} $ is  defined
for some  $n\geq 2$.  We define  $\mu
^{(n+1)}$  as the
composite  $
\Omega B (A^{\otimes n+1}) \stackrel{\theta _{_{A, A^{\otimes n}}}}{\to}
 \Omega B \left( A \otimes \Omega B (A^{\otimes n})\right)
\stackrel{ \Omega B(id \otimes  \mu ^{(n)})}{\to}
 \Omega B \left( A \otimes \Omega B (A)\right)
\stackrel{\Omega B(id \otimes \alpha_A) }{\to}
 \Omega B \left( A \otimes A\right)
\stackrel{ \mu }{\to}
\Omega B A $.
If we assume moreover that
 $\alpha _{A} \circ \mu ^{(n)} \circ i_{A^{\otimes n}} = m_A^{( n)}$
 then the identities ($\star$) and
$\alpha_{-}
\circ i_{-} = id _{-}$ imply that $\alpha _{A}
\circ \mu ^{(n+1)} \circ i_{A^{\otimes n+1}} = m _A^{(n+1)}$. The last
statement is Proposition  4.5  in \cite{[Mu]}.
\vspace{3mm }

\vspace{3mm }

\noindent{\bf A.4 Lemma.} {\it Let $(W, \psi _W)$ be the standard resolution.

a) If $(A,d)$ is a differential graded algebra then    $ \mbox{{\rm
Hom}}(W, A)\stackrel{ev_0}{\to} A \in {\bf TEX}_A$.

b) If $ X \stackrel{\alpha}{\to}  A \in {\bf TEX}_A $ then $ \mbox{\rm
Hom}(W,X) \stackrel{\mbox{\rm Hom}(W, \alpha ) }{\to}
\mbox{\rm Hom}(W, A) \in {\bf TEX}_{\mbox{\tiny \rm Hom} (W,A)}$.

c) If $A$ is a $shc$ algebra then

\hspace{5mm} (i) $\mbox{\rm Hom}(W,A)$ is a $shc$ algebra,

\hspace{5mm}  (ii) $ev_0$ is a $shc$ homomorphism,

\hspace{5mm}  (iii) Let  $A'$ be  a $shc$ algebra.  A homomorphism of
differential graded algebras\\
 $f: A'\rightarrow \mbox{\rm Hom}(W,
A)$ is a $shc$ homomorphism if and only if
$ev_0\circ f$  is a  $shc$ homomorphism.}

\vspace{2mm}

\noindent{ \bf Proof.} a) Let us denote $[g_0,g_1,...,g_n] \in G{^\times
n+1}$ a generator on $W_n$ and  recall the linear
homotopy
$k: W\rightarrow W$ defined by
$k([g_0,...,g_n]) = \sum_{i=0}^n (-1)^i[g_0,..,g_{i-1}, 1, g_i, ..., g_n]$.
Obviously,
 $k\circ \partial_W +\partial_W\circ k = id_W - \eta_W\circ \epsilon_W$ and
$k\circ k =0$.
We define $h: \mbox{Hom}(W, A) \rightarrow \mbox{Hom}(W, A)$ by $h(f) =
k\circ f$.  If $ \rho :  A \rightarrow \mbox{Hom}(W,
A)$ is defined by
$\rho (a)(w)= a\eta_A\epsilon_W(w)$, $ a \in A$,  $w\in W$,   a
straightforward computation shows that
 $h\circ D + D\circ h = id_{\mbox{\tiny Hom}(W,A)} - \rho \circ ev_0$ and
$h\circ h =0$.   One also  easily checks that
$ev_0\circ\rho = id_A$,
$ ev_0\circ h =0$, $h\circ\rho = 0$,
$\rho\circ\eta_A =\eta_{\mbox{\tiny Hom}(W,A)}$, $\epsilon_{\mbox{\tiny
Hom}(W,A)}\circ \rho =
\epsilon_A$,  where $ \eta_{\mbox{\tiny Hom}(W,A)} = \eta_A\circ
\epsilon_W$ and
$\epsilon_{\mbox{\tiny Hom}(W,A)}(f) =
\epsilon_A\circ f\circ \eta_W(1).$

\vspace{2mm}
b) The proof is similar to the end of the  proof of part a).

\vspace{2mm}
c) (i) Consider the differential linear map
$
\psi_A:  \mbox{Hom} (W, A)\otimes \mbox{Hom}(W, A) \rightarrow
\mbox{Hom}(W, {A\otimes A})
$
defined by
$\psi_A (f\otimes g) = (f\otimes g) \circ \psi_W$. Then,   $\mbox{Hom}(W,
m_A)\circ \psi_A =
\cup_{\mbox{\tiny Hom}(W, A)}$ and the  following equalities prove that
$\psi_A$ is a homomorphism of differential graded algebras.

\centerline{$
\begin{array}{rcl}
\psi_A ((f_1\otimes g_1)\cup (f_2\otimes g_2))& = &
(-1)^{\vert f_2\vert \vert g_1\vert}\psi_A  ((f_1\cup f_2) \otimes (g_1\cup
g_2))
\\
  &  =  & (-1)^{\vert f_2\vert \vert g_1\vert} ((f_1\cup f_2)\otimes (g_1\cup
g_2))\circ \psi_W
\\
& = &(-1)^{\vert f_2\vert \vert g_1\vert} (m_A\circ (f_1\otimes  f_2)\circ
\psi_W)\otimes  m_A\circ ((g_1\otimes  g_2)\circ \psi_W) \circ \psi_W
\\
 & = & m_{A\otimes A} \circ ((f_1\otimes  g_1)\otimes  (f_2\otimes  g_2)
)\circ  \psi_W\circ \psi_W
\\
  &= & m_{A\otimes A} \circ (\psi_A (f_1\otimes  g_1)\otimes \psi_A (
f_2\otimes   g_2))\circ \psi_W
\\
  &= &\psi_A (f_1\otimes g_1)\cup \psi_A(f_2\otimes g_2).
\end{array}
$}

\noindent
 We define $\mu _{\mbox{\tiny Hom}(W, A)}: \Omega B \left( \mbox{Hom}(W ,
A)^{\otimes 2}\right) \to \Omega B
\mbox{Hom}(W , A)$    as  the composite
 $  \Omega B\mbox{Hom}(W, \alpha_A\circ  \mu_A)\circ
\theta_{\mbox{\tiny Hom}(W, A\otimes A)}\circ\Omega B\psi$ where
$\theta_{\mbox{\tiny Hom}(W, A\otimes A)}$ denotes the
unique homomorphism of differential graded algebras such that
$\theta_{\mbox{\tiny Hom}(W, A\otimes A)}\circ
\mbox{Hom}(W, \alpha_A\otimes \alpha _A) \circ \alpha _{\mbox{\tiny Hom}(W,
\Omega B A \otimes \Omega B A)} =
\alpha _{\mbox{\tiny Hom}(W,  A \otimes  A)}$. Existence and  unicity of
$\theta_{\mbox{\tiny Hom}(W, A\otimes A)}$ is a
direct consequence of F3, F1, F2 and part b). It  turns out that $(\mbox{
Hom}(W, A),\mu _{\mbox{\tiny Hom}(W,
A)})  $ is a  $shc$ algebra.

\vspace{2mm}
c)-(ii). Consider for any $A$ in {\bf DA} the commutative diagram

\centerline{$
\begin{array}{ccccc}
\Omega BA &\stackrel
{\Omega B(ev_0)}{\longleftarrow} & \Omega B\mbox{Hom}(W,A)&\stackrel{
\alpha_{\Omega B\mbox{\tiny Hom}(W,A)}}{\longleftarrow}
& (\Omega B)^2\mbox{Hom}(W,A)
\\
 &&\hspace{-15mm} \alpha_{\mbox{\tiny Hom}(W,A)}\downarrow && \hspace{15mm}
\downarrow \Omega B\alpha_{\mbox{\tiny Hom}(W,A)}
\\
&\hspace{-10mm}\alpha_A \searrow &\mbox{Hom}(W,A)&\stackrel
{\alpha_{\mbox{\tiny Hom}(W,A)}}{\longleftarrow} &\Omega
B\mbox{Hom}(W,A)B\\ &&ev_0\downarrow&&\downarrow  \Omega Bev_0
\\
&&A&\stackrel
{\alpha_A}{\longleftarrow} &\Omega BA \,.
\end{array}
$}

\noindent
>From  F3,  we deduce the existence of  $\theta'_{\Omega BA} : \Omega B  A
\rightarrow \Omega B \mbox{Hom}(W,A)$ such that
$ev_0\circ \alpha_{\mbox{\tiny Hom}(W,A)}\circ \theta'_{\Omega BA}=
\alpha_A$ and
$\theta'_{\Omega BA}\circ \Omega B\mbox{Hom}(W,ev_0)\simeq_{{\bf DA}}
\alpha_{\Omega B\mbox{Hom}(W,A)}$. By  F7,
$\Omega B\alpha_{\mbox{\tiny Hom}(W,A)}\simeq_{{\bf DA}}\alpha_{\Omega
B\mbox{\tiny Hom}(W,A)}$ and by  F5, there is  a
section for
$\alpha_{\Omega B\mbox{\tiny Hom}(W,A)}$ so that
$ (\star_A) \quad \quad \theta'_{\Omega BA}\circ\Omega Bev_0
\simeq_{{\bf DA}} id_{\Omega B\mbox{Hom}(W,A)}\,.
$
Moreover, from the unicity property  in  F3, we know that
$\Omega Bev_0\circ \theta'_{\Omega BA} = id_{\Omega BA}.$  In  the
following diagram,  the left hand square is commutative,

\centerline{$
\begin{array}{ccccc}
\Omega BA &\stackrel
{\Omega B(ev_0)}{\longleftarrow} & \Omega B\mbox{Hom}(W,A)&\stackrel
{\theta'_{\Omega BA}}{\longleftarrow} & \Omega BA\\
\mu_A\uparrow& &\hspace{20mm} \uparrow \Omega B\mbox{Hom}(W,\alpha_A
\circ\mu_A )&&\uparrow\mu_A\\
\Omega B\left(A^{\otimes 2}\right)&\stackrel
{\alpha_{\Omega B\left(A^{\otimes 2}\right)}\circ\Omega B(ev_0)
}{\longleftarrow} &
 \Omega B\mbox{Hom}\left(W,\Omega B\left(A^{\otimes 2}\right)\right)&\stackrel
{\theta'_{\Omega B\left(A^{\otimes 2}\right)}}{\longleftarrow} & \Omega
B\left(A^{\otimes 2}\right)\,.
\end{array}
$}

\noindent
Since  the composite of  horizontal maps  are the
identity maps, relations
$(\star_A)$ and $(\star_{A\otimes A })$ imply that the right hand square
commutes up to
homotopy in {\bf DA}. Therefore, in the following diagram the upper square
commutes , up to homotopy in {\bf DA}.

\centerline{$
\begin{array}{ccccc}
\Omega B\mbox{Hom}(W, A) &\stackrel {\Omega B(ev_0)}{\longrightarrow} &
\Omega BA&\stackrel{\theta'_{\Omega
BA}}{\longrightarrow} & \Omega B\mbox{Hom}(W,A)
\\
|| & &&&  \uparrow  \theta'_{\Omega BA} \circ \mu_A
\\
\Omega B\mbox{Hom}(W, A) &\stackrel{\Omega B\mbox{\tiny Hom}(W,\alpha_A
\circ\mu_A )}{\longrightarrow}
 &\Omega B\mbox{Hom}(W, \Omega B\left( A^{\otimes 2}\right) )&
\stackrel{\theta'_{\Omega B(A^{\otimes 2})}}{\longleftarrow} &\Omega
B(A^{\otimes 2})
\\
\hspace{-11mm}\mu_{\mbox{\tiny Hom}(W,A)} \uparrow && \hspace{-20mm} \Omega
B\mbox{Hom}(W,\Omega B(A^{\otimes 2}))\downarrow
&& ||
\\
\Omega B(\mbox{Hom}(W,A)^{\otimes 2})&\stackrel
{\Omega B\psi_A }{\longrightarrow} &\Omega B\mbox{Hom}(W, A^{\otimes 2})&
\stackrel{\theta'_{\Omega B(A^{\otimes 2})}}{\longleftarrow}&\Omega
B(A^{\otimes 2})
\\
\hspace{-20mm}\Omega B(ev_0\otimes ev_0)\downarrow & & \Omega
B(ev_0)\downarrow && ||
\\
\Omega B(A^{\otimes 2})&=& \Omega B(A^{\otimes 2})&=&\Omega B(A^{\otimes 2})
\end{array}
$}

\noindent
It is not difficult, if tedious, to check that the other cells in the
diagram commutes also  up to homotopy in {\bf DA}. This
shows that  $\Omega B(ev_0) \circ \mu_{\mbox{\tiny Hom}(W,A )} \simeq _{\bf
DA} \mu _A \circ \Omega
B(ev_0\otimes ev_0)$, i.e.  $ev_0$ is a strict $shc$ map.

\vspace{2mm}
c)-(iii)  Let $f : A' \to \mbox{Hom}(W,A) $ be a homomorphism of
differential graded algebras such that
$\mu_A \circ \Omega (ev_0\otimes ev_0) \circ \Omega B (f\otimes f) \simeq
_{\bf DA} (\Omega B (ev_0 \circ f)) \circ
\mu_{A'}$. Thus in the following diagram the bigger square commutes, up to
homotopy in ${\bf DA}$.

\centerline{$
\begin{array}{ccccc}
\Omega BA' &\stackrel
{\Omega Bf}{\longrightarrow} & \Omega B\mbox{Hom}(W,A)&\stackrel
{\Omega Bev_0}{\longrightarrow} & \Omega BA\\
\mu_{A'}\uparrow& & \hspace{-6mm}\mu_{\mbox{\tiny Hom}(W,A)}\uparrow  &
&\uparrow\mu_A\\
\Omega B(A'\otimes A')&\stackrel
{\Omega B(f\otimes f)}{\longrightarrow} &
 \Omega B(\mbox{Hom}(W,A)\otimes \mbox{Hom}(W,A))&\stackrel
{\Omega B(ev_0\otimes ev_0)}{\longrightarrow} & \Omega B(A\otimes A)\,.
\end{array}
$}

\noindent
The right hand square commutes, up to homotopy in ${\bf DA}$, by part
c-(ii) and  thus so does the left hand square.  This
proves that   $f$ is a strict $shc$ homomorphism. Conversely, if $f$ is a
$shc$ map, by part c-(ii), the composite $f\circ
ev_0$ is also a
 strict $shc$ homomorphism.

\vspace{3mm }

\noindent{\bf A.5 Lemma.} {\it  If $(A, d_A,\mu_A, \tilde\kappa_A)$ is a
$\pi$-$shc$ algebra then,

a) $\tilde \kappa $ is a strict homomorphism of $shc$  algebras,

b) $H^*\phi_A: H{\frak C}(\mbox{\rm Hom}(W, A))\rightarrow H(\mbox{\rm
Hom}(W,{\frak C}A))$
preserves the natural multiplications.}

 \vspace{2mm }

\noindent{\bf Proof.} a) By lemma A.3,  $\mu^{(p)}$ is a strict homomorphism
of $shc$ algebras.  Since $ev_0\circ \tilde \kappa_A \simeq_{{\bf DA}}
\alpha_A\circ\mu_A^{(p)}$,
from
 Lemma A.4-c)-(iii), we deduce that $\tilde \kappa_A $ is a strict $shc$
homomorphism.

b)  We obtain   the two  commutative diagrams (Diagrams A and B):

\centerline{$
\begin{array}{ccccccccccccc}
\displaystyle 
\left( \frak{C}\mbox{Hom}(W,A)\right)^{ \otimes 2}
 & \stackrel{\phi _A ^{\otimes 2} }{\longrightarrow}
& \mbox{Hom}(W,{\frak C} A) ^{\otimes 2}
 &\stackrel{ \psi _{{\frak C} A }^{ \otimes 2} }{\longrightarrow}
& \mbox{Hom}\left( W,({\frak C}A)^ { \otimes 2} \right)
 \\
\displaystyle 
 \hspace{5mm} \downarrow sh  &&&& \hspace{-25mm}\mbox{Hom}(W,sh) \downarrow
\\
\displaystyle 
 \frak{C}\left(\mbox{Hom}(W,A)^{ \otimes 2} \right)
 &\stackrel{{\frak C} \psi _{\mbox{\tiny  Hom}(W,A) }} {\longrightarrow}
& {\frak C}\mbox{Hom}(W,A ^{\otimes  2}
 &\stackrel{ \phi _{A^{ \otimes 2}}   }{\longrightarrow}
&  \mbox{Hom}(W,{\frak C} (A^{\otimes 2}))
\\
\displaystyle 
\hspace{-20mm} {\frak C}\alpha _{\mbox{\tiny  Hom}(W,A)  ^{\otimes 2}}
\uparrow
&& \hspace{15mm}\uparrow
{\frak C } \mbox{Hom}(W,\alpha _{A^{\otimes 2}})
&&
\\
\displaystyle 
\frak{C} \Omega B \left( \mbox{Hom}(W,A)^{\otimes 2}  \right)
 &\stackrel{{\frak C} \psi _{ \Omega B(\mbox{\tiny  Hom}(W,A)^{ \otimes
2}}} {\longrightarrow}
& {\frak C}\Omega B \mbox{Hom}(W,A^{ \otimes 2} )
&
\,,
\end{array}
$}

\noindent and

\centerline{$
\begin{array}{ccccccccccccc}
\displaystyle 
 {\frak C}\mbox{Hom}(W,A^{\otimes  2})
 &\stackrel{ \phi _{A ^{\otimes 2}}   }{\longrightarrow}
&  \mbox{Hom}(W,{\frak C} (A^{\otimes 2}))
\\
\displaystyle 
\hspace{-20mm}{\frak C } \mbox{Hom}(W,\alpha _{A^{\otimes 2}}) \uparrow
&& \hspace{20mm}\uparrow \mbox{Hom}(W,{\frak C} \alpha _{A^{\otimes 2}})
\\
\displaystyle 
{\frak C}\mbox{Hom}(W, \Omega B (A^{\otimes 2}))
&\stackrel{ \phi _{\Omega B(A ^{\otimes 2})}   }{\longrightarrow}
&  \mbox{Hom}(W,{\frak C}\Omega B (A^{\otimes 2}) \,.
\end{array}
$}

\noindent
From, F3, Lemma A.4-b) and  F8 we deduce the existence of the homomorphisms
of differential graded algebras $\theta
'_A$,
$\theta '_{A^{\otimes 2}} $, $\theta ''_{A^{\otimes 2}} $
 such that : $\mbox{Hom}(W, \alpha _A) \circ \theta ' _A = \alpha
_{\mbox{\tiny Hom} (W,A)} $, $\mbox{Hom}(W, \alpha
_{A^{\otimes 2}}) \circ \alpha _{\mbox{\tiny Hom} (W,\Omega B(A^{\otimes
2}))} \circ  \theta '' _{A^{\otimes 2}}  = \alpha
_{\mbox{\tiny Hom} (W,A^{\otimes 2})} $ and
 $\mbox{Hom}(W, \alpha _{A^{\otimes 2} }) \circ \theta ' _{A^{\otimes 2}}
= \alpha _{\mbox{\tiny Hom} (W,A^{\otimes 2})} $.
The desired commutativity of the  diagrams C, and D below  is clear from
the construction of $ \mu _{\mbox{\tiny Hom}( W,A)
}$ (Proof of lemma A.4-c)-(i)):

\centerline{$
\begin{array}{cccccccc}
\displaystyle 
\frak{C} \Omega B \left( \mbox{Hom}(W,A)^{ \otimes 2} \right)
 &\stackrel{{\frak C} \psi _{ \Omega B(\mbox{\tiny  Hom}(W,A)^{ \otimes 2})
}} {\longrightarrow}
& {\frak C}\Omega B \mbox{Hom}(W,A^{ \otimes 2} )
\stackrel{{\frak C} \theta '' _{A^{\otimes 2} }}{\longrightarrow} &{\frak
C}\Omega B \mbox{Hom}(W,\Omega B (A ^{\otimes 2}))
\\
\displaystyle 
\hspace{10mm}{ \frak C}\mu _{\mbox{\tiny  Hom}(W,A)} \downarrow
&&
&\hspace{-10mm} {\frak C} \alpha _{\mbox{\tiny  Hom}(W,\Omega B A)}
\circ \downarrow { \frak C} \Omega B \mbox{Hom}(W,{\frak C}\mu _A)
\\
 \displaystyle 
{\frak C}\Omega B \mbox{Hom}(W,A)
&& \stackrel{{\frak C} \theta ' _A}{ \longrightarrow }
&  {\frak C} \mbox{Hom}(W,\Omega BA)\,,
\end{array}
$}

\noindent
and

\centerline{$
\begin{array}{cclcccccccccc}
\displaystyle 
 {\frak C}\Omega B \mbox{Hom}(W,A ^{\otimes 2})
&\stackrel{ {\frak C}\theta ' _{\Omega B (A ^{\otimes 2})}   }{\longrightarrow}
&{\frak C} \mbox{Hom}(W, \Omega B (A^{\otimes 2}))
&\hspace{-5mm}\stackrel{ \phi _{\Omega B(A ^{\otimes 2})}   }{\longrightarrow}
& \hspace{-5mm} \mbox{Hom}(W,{\frak C}\Omega B (A^{\otimes 2})
\\
\displaystyle 
\hspace{-10mm}   {\frak C}\theta '' _{ A ^{\otimes 2}} \downarrow
& \hspace{-22mm}   {\frak C} \alpha_{ \mbox{\tiny Hom}(W, \Omega B
(A^{\otimes 2}))}\nearrow
&&&
\\
\displaystyle 
  {\frak C}\Omega B \mbox{Hom}(W,\Omega B (A^{ \otimes 2}))
&&&& \hspace{10mm} \downarrow  \mbox{Hom}(W,{\frak C}\mu _A)
\\
\displaystyle 
\hspace{5mm}  {\frak C} \alpha _{\mbox{\tiny  Hom}(W,\Omega B A)}
\circ \downarrow { \frak C} \Omega B \mbox{Hom}(W,{\frak C}\mu _A)
&&&&
\\
\displaystyle 
 {\frak C} \mbox{Hom}(W,\Omega BA)
&& \stackrel{ \phi _{\Omega B A}}{ \longrightarrow }
&&  \mbox{Hom}(W,{\frak C}\Omega BA)\,.
\end{array}
$}

\noindent
Now we choose linear sections of ${\frak C} \alpha _{\mbox{\tiny Hom}(W, A
^{\otimes 2})} $ (resp. of
$ \mbox{ Hom}(W,{\frak C} \alpha _{ A ^{\otimes 2}) }$) and  we define a
product $m_{ {\frak C}\mbox{\tiny Hom}(W,
A)}$  (resp. a cup product $ m_{\mbox{\tiny Hom}(W, {\frak C}A)})$ on
${\frak C}\mbox{Hom}(W,
A)$ (resp. on
$\mbox{Hom}(W, {\frak C}A)$. Then gluing together diagrams A,B,C, and D we
deduce
that  the following diagram

\centerline{$
\begin{array}{cccc}
\left( {\frak C}\mbox{Hom}(W, A) \right) ^{\otimes 2}  &
\stackrel{\phi_A^{\otimes 2}}{\longrightarrow}
& \left( \mbox{Hom}(W, {\frak C}A) \right) ^{\otimes 2}
\\
\hspace{-15mm}  m_{ {\frak C}\mbox{\tiny Hom}(W, A)} \downarrow & &
\hspace{15mm} \downarrow m_{\mbox{\tiny Hom}(W, {\frak C}A)})
\\
 {\frak C}\mbox{Hom}(W, A)   & \stackrel{\phi_A}{\longrightarrow}
&  \mbox{Hom}(W, {\frak C}A)
\end{array}
$}

\noindent
commutes up to a linear  homotopy.

\vspace{3mm}

\noindent{\bf A.6 Lemma.}  {\it   Let $(TV, d_V) $  be  a differential
graded algebra   and
assume that a finite group
$\pi$ acts freely on  $V$.  Let
$A$ be a $\pi$-differential graded algebra  and $f, g\in \pi$-${\bf DA}(TV,
A)$. If
$f\simeq _{\pi\! -\!{\bf DA}} g$  then
 ${\frak C}f\simeq _{\pi }{\frak C}g$.}
 \vspace{2mm }

\noindent{\bf Proof.} The closed model category
framework  provides a convenient language in which we prove the lemma.
For our purpose we define the cylinder
$I(TV, d)$ on
$(TV, d)$:

\centerline{$
\begin{array}{ccccc}
(TV, d) &\searrow \partial_0&&
 &
 \\
 &&I(TV, d)&
 \stackrel {p}{\longrightarrow}&
 (TV, d)\\
(TV, d) &\nearrow \partial_1&&
 &
\end{array}
$}

\noindent
with $ I(TV, d) := (T(V_0\oplus  V_1\oplus sV), D)$, $\partial_0V=V_0$,
$\partial_1V=V_1$, $pv_0=pv_1=v$, $psv=0$, $D=d$ on $V_0$ and on $V_1$,
$Dsv=Sdv$ where
 $S$ is the unique
$(\partial_0, \partial_1)$-derivation
$S: TV\rightarrow T(V_0\oplus  V_1\oplus sV)$ extending the graded isomorphism
$s:V\to sV$. The free $\pi$-action on $TV$ naturally extends to a free
$\pi$-action on
$I(TV, d)$ so that $I(TV, d)$ is a $\pi$-algebra and the maps
 maps $p, \partial_0$ and   $\partial_1$ are $\pi$-equivariant
 quasi-isomorphisms. By definition $f\simeq _{C}g$ (resp. $f\simeq
_{\pi\mbox{-}C}g$) if there exists $H\in  {\bf DA}( I(TV,
d), A)$ (resp.
$H\in  \pi$-${\bf DA}( I(TV, d), A)$)
such that $H\partial_0 = f$ and $H\partial_1 = g$. It is straightforward to
check that
$f\simeq _{C}g$ if and only if $f\simeq _{{\bf DA}}g$ , \cite{[FHT2]},
(resp. $f\simeq _{\pi\mbox{-}C}g$ if and only if
$f\simeq _{\pi\mbox{-}{\bf DA}}g$). Consider  the commutative diagram

\centerline{$
\begin{array}{ccc}
{\frak C}TV\oplus {\frak C}TV &
 \stackrel {{\frak C}id\oplus {\frak C}id}{\longrightarrow}&
 {\frak C}TV\\
 {\frak C}\partial_0\oplus {\frak C}\partial_1\downarrow&&{\frak C}p\uparrow\\
 {\frak C}ITV&
 \stackrel {id }{\longrightarrow} &{\frak C}ITV.\\
\end{array}
$}

\noindent
Now ${\frak C}\partial_0\oplus {\frak C}\partial_1$ is injective and its
cokernel is  a projective
$\pi$-module.  Hence, ${\frak C}\partial_0\oplus {\frak C}\partial_1$  is a
cofibration and since
${\frak C}p$ is a weak equivalence, the complex   ${\frak C}ITV$ is a
cylinder object in the closed model category $\pi$-${\bf
DM}$.  Let
$H\in  \pi$-${\bf DA}( I(TV, d), A)$ be a homotopy between $f$ and $g$
then, ${\frak C}f= {\frak C}H
\circ {\frak C} \partial _0$, ${\frak C}g= {\frak C}H \circ {\frak
C} \partial _1$ and thus  ${\frak C}f \simeq_{\pi} {\frak C}g$.

\vspace{ 1cm}

\centerline{ \large Appendix B - Equivariant acyclic model theorem. }

\vspace{3mm }
The proof of Theorem B  relies  heavily on the  $\pi
$-equivariant acyclic model theorem  for cochain
 functors.  We  state and prove this theorem.

\vspace{3mm}

\noindent{\bf B.1} Let $R$ be a (ungraded) commutative algebra over the field
${\Bbb F}_p$ and let ${\bf C}$ be a category with
models
${\cal M}$.  Consider a contravariant functor
$
F :\mbox{\bf C} \to \mbox{ \bf Coch}_R\,, A \mapsto FA = \{F^i A\}_{i \in
{\Bbb
Z}}
$  with values in the category  of  $R$-cochain
complexes ($R$ acting on the left). See \cite{[BG]}, \cite {[Mw]} or
\cite{[BT]} for the definitions of:
 $F$  admits a unit,  $F$ is acyclic on the models and  $F$ is
corepresentable for
the models ${\cal M}$.   The singular cochains  functor
$ C^n : \mbox{\bf Top} \to \mbox{\bf Mod }_R $
is corepresentable on the models $\{\triangle ^k\}_{k\geq 0}$.  Observing
that a retract  of
corepresentable functor is corepresentable,
we deduce  that the functor $X \mapsto N^*X$  is corepresentable on the models
$\{\triangle ^k\}_{k\geq 0}$. For further use it is interesting to  remark
here
that is $F$ is corepresentable on a familly of
models  and if $W$ is an $R$-free graded module of finite type then so is the
functor $A \mapsto \mbox{Hom}(W,FA)$. Indeed, in this
case  each $ \mbox{Hom}(W_n ,FA)$  is a finite sum on copies of $FA$ and
one concludes
using the fact that $\widehat {\prod_{\alpha }
F_\alpha }= {\prod_{\alpha } \hat F_\alpha }$ (see \cite{[BT]}-$I\!I$ for
more details).

\vspace{3mm}

\noindent{\bf B.2} We embed the category ${\bf Coch}_R$ into the category
${\bf
Mod}$ of graded (without differential) ${\Bbb F}_p$-modules. Then,
considering  the two  contravariant functors
$F\,, G  :\mbox{\bf C} \to \mbox{ \bf Coch}_R
$
as functors  with values in $
\mbox{ \bf Mod}_R $,  we
denote by
$ \mbox{Hom}^i( F,G)\,, \quad  i \in {\Bbb Z}  $
 the sets of natural transformations of degree $i$.  A
differential
$D : \mbox{Hom}_R^i ( F, G)  \to \mbox{Hom}^{i+1} ( F, G)
$
is defined by: $(DT)_A= d_{G(A)} T_A - (-1)^iT_A d_{F(A)}$.  For our purpose
  $R= {\Bbb F}_p[\pi] = {\Bbb F}_p[ \tau ]/( \tau
^p-1)$ then, for any object $A$ of
${\bf C}$,
$
\mbox{Hom}^i ( FA, GA)
$ is a
$\pi$-module by the rule
$
 \sigma  \in \pi \,,  T \in \mbox{Hom}_R^i ( F, G) \,, A, \in {\bf C} \,, x
\in
FA \,,\quad  (\sigma T)_A (x) = \sigma
T_A(\sigma ^{-1} x)\,.
$
The ${\Bbb F}_p $-module $\mbox{Hom}_\pi  ( FA, GA)$ of $\pi $-linear
transformations is the  fixed point set of
$\mbox{Hom } ( FA, GA)$ under this action. The next theorem hase been prove
in \cite{[BG]}.

\vskip 3mm
\noindent {\bf B.3  Theorem.}
{\it Let ${\bf  C}$ be a category with models ${\cal M}_R$ and
$F, G: {\bf  C}\rightarrow {\bf Coch}_R $
two contravariant functors with units such that  for any $A\in {\bf C}$,
$F^i A = 0 = G^i A$ if $i<0$,
 $F$ is acyclic on the models and  $G$ is corepresentable on the models. Then
$
H^0 \left(\mbox{\rm Hom}_\pi ( F, G)\right)=
R \mbox{ and  } H^0 \left(\mbox{\rm Hom}_\pi ( F, G)\right)=
0 \mbox{ if } i \neq 0\,.
$
}

\vskip 3mm
\noindent {\bf B.4 Corollary.} {\it If $ F: \mbox{\bf Top} \to \mbox{\bf
Coch}_\pi $ is any functor such that:

a) $\eta :  {\Bbb F}_p \to F $ is a unit,

b) for any model $\triangle ^n \,, n \geq 0$ there exists a map $\epsilon _n :
F(\triangle ^n ) \to {\Bbb F}_p $ such that
$\epsilon _n \circ \eta _n = id_{{\Bbb F}_p} $ and $\eta_n \circ \epsilon
_n \simeq
id _{F\triangle ^n )}$ then,
$
H^i( \mbox{\rm Hom}_\pi (F, \mbox{\rm Hom}(W,N^*)) = \left\{
\begin{array}{lr}
{\Bbb F}_p[\pi]  &i=0 \\
0  &i >0
\end{array}\,.
\right.
$
}

\vspace{2mm}
\noindent {\bf Proof}    Recall
the (right)
$\pi$-free acyclic complex
$W$ with left  $\pi $-action defined by $ \sigma w = w \sigma ^{-1} \,, \quad
\sigma \in \pi$ and $ w \in W$ and the canonical
isomorphism
$
{\mbox{Hom}}_\pi \left( A,  {\mbox{Hom}}(W, B))\right)  \cong
{\mbox{Hom}}_\pi \left(W \otimes   A,
B)\right) \,,
$
where $(B, d_B)$ is a (left)   $\pi$-complex. We precise  that $\pi$ acts
tivially on $N^* X$ for any space $X$.  The functor
$
G=\mbox{Hom}(W, N^*) : \mbox{\bf Top} \to \mbox{ \bf Coch}_{{\Bbb F}_p
[\pi]}\,,$
admits a unit $\eta $ defined as follows: for any space $X$
$ \eta _X ( \sigma ) = \sigma \circ e^*_0 \,,
$
where $e^*_0 \in \mbox{Hom } (W, N^*X) $ denotes  the dual of $e_0 \in
W_0$.  As  remarked  in B.1,  $G=\mbox{Hom}(W, N^*)$
is  corepresentable   on the models $\triangle ^n$.

 The natural transformation $\eta $ induces a natural
transformation, also denoted $\eta : R= {\Bbb F}_p[\pi] \to F \otimes R$.
Moreover, for any $n\geq 0$ there exists $\epsilon
_n  = F(\triangle ^n ) \otimes R \to R $ such
that  $
\eta _{\triangle ^n} \circ \epsilon _{\triangle ^n} = id_{\triangle ^n}$  and $
\epsilon _{\triangle ^n} \circ \eta _{\triangle ^n} \simeq id_{\triangle ^n}$.
Since  $\pi $ acts diagonally on  $F \otimes R$, the map $ \varphi \mapsto
(\hat \varphi : x\otimes y \mapsto \varphi (x)y)
$  is an isomorphism of chain complexes  $
\mbox{Hom} _{\pi}(F, G) \cong \mbox{Hom} _{\pi} ( F \otimes R, G)\,, $
for
any functor  $G : \mbox{\bf Top} \to \mbox{\bf Coch}_{\pi} $. We apply theorem
B.3 to end the proof.


\vspace{1 cm}

\begin{thebibliography}{99}




 \bibitem{[B]} D.J. Benson, {\it Representations  and Cohomology I$\!$I},
Cambridge studies in advanced mathematics {\bf 31} - Cambridge University
Press - 1991.


 \bibitem{[BT]} N. Bitjong  and J-C. Thomas, {\it On the cohomology algebra of
free loop spaces},  Topology (2001) to appear.

 \bibitem{[BG]} A.K. Bousfield and V.K.A.M. Gugenheim, {\it On PL de
Rham theory and rational homotopy type}, Memoirs of the Amer.
       Math. Soc. {\bf 179} 1976.

\bibitem{[C-E]} H. Cartan and S. Eilenberg, {\it Homological algebra},
Princeton mathematical series {\bf 19} -
Princeton University Press - 1956.







 \bibitem{[Eil]} S. Eilenberg, {\it Topological methods in abstract algebra},
Bull. Amer. Math. Soc. {\bf 55} (1949) 3-27.

 \bibitem{[FHT]} Y. F\'elix, S. Halperin and J-C. Thomas, {\it Adams'
       cobar equivalence}, Trans. of the Amer. Math. Soc. {\bf 329} (1992)
531-549.


\bibitem{[FHT2]} Y. F\'elix, S. Halperin and J-C. Thomas,
{\it Differential graded algebras in topology},  Handbook of
Algebraic Topology, I. James,  editor,  Chapter 16 - North-Holland - 1995.

\bibitem{[HL]} S. Halperin and J. M. Lemaire,
{\it Notions  of category in differential algebras}, Algebraic Topology:
Rational Homotopy,  Springer Lecture Notes in Math. {\bf 1318} (1988) 138-153.

\bibitem{[HV]} S. Halperin and M. Vigué-Poirrier, {\it The homology of a free
loop space}, Pacific J. of Math. {\bf 147} (1991) 311-324.



\bibitem{[Jo]} J. D.S. Jones, {\it Cyclic homology and equivariant homology},
Invent. Math. {\bf 87} (1987) 403-423.




  \bibitem{[L]} J-L. Loday, {\it Cyclic homology}, Grundleren der
mathematischen
       Wissenschaften {\bf 301} - Springer-Verlag - 1991.




\bibitem{[Mw]} M. Majewski, {\it Rational homotopical models and uniqueness},
Memoirs of the Amer. Math. Soc. {\bf 143/682}, 2000.


 \bibitem{[Ma]} J.P. May, {\it A general algebraic approach to Steenrod
operations}, Springer Lecture Notes in Math {\bf 168}  (1970) 153-231.

   \bibitem{[ML]} S. Mac Lane, {\it Homology}, Grundleren der mathematischen
       Wissenschaften  {\bf 114} Springer-Verlag - 1963.

   \bibitem{[Mu]} H. J. Munkholm, {\it The Eilenberg-Moore spectral
sequence and
       strongly homotopy multiplicative maps}, J. of Pure and Appl. Alg.
{\bf 5} (1974) 1-50.




 \end{thebibliography}
\end{document}